\documentclass[preprint,12pt]{elsarticle}
\usepackage{amsfonts, amsmath, amsthm, hyperref, caption, subcaption,longtable}
\usepackage[bbgreekl]{mathbbol}
\usepackage{stmaryrd}
\usepackage{tikz}
\usepackage{graphicx} 
\usepackage{url}
\usepackage{hyperref}
\usepackage[numbers]{natbib}
\usepackage{mathrsfs}
\usepackage[ruled,linesnumbered,noend]{algorithm2e}
\usepackage{bm} 

\usepackage{cleveref}
\crefformat{table}{#2Tab.~#1#3}
\Crefformat{table}{#2Table~#1#3}

\crefformat{figure}{#2Fig.~#1#3}
\Crefformat{figure}{#2Figure~#1#3}

\crefformat{equation}{#2Eq.~(#1)#3}
\Crefformat{equation}{#2Equation~(#1)#3}
\crefmultiformat{equation}{Eqs.~(#2#1#3)}%
{ and~(#2#1#3)}{, (#2#1#3)}{ and~(#2#1#3)}

\Crefmultiformat{equation}{Equations~(#2#1#3)}%
{ and~(#2#1#3)}{, (#2#1#3)}{ and~(#2#1#3)}

\crefrangeformat{equation}{#3Eqs.~(#1)#4--#5(#2)#6}
\Crefrangeformat{equation}{#3Equations~(#1)#4--#5(#2)#6}

\crefformat{section}{#2Section~#1#3}
\Crefformat{section}{#2Section~#1#3}

\crefformat{subsection}{#2Subsection~#1#3}
\Crefformat{subsection}{#2Subsection~#1#3}

\crefformat{chapter}{#2Chapter~#1#3}
\Crefformat{chapter}{#2Chapter~#1#3}

\Crefformat{appendix}{#2#1#3}
\crefformat{appendix}{#2#1#3}

\crefformat{algorithm}{#2Alg.~#1#3}
\Crefformat{algorithm}{#2Algorithm~#1#3}

\usepackage{lineno}

\journal{Computers \& Structures}

\def\rev#1{{#1}}

\numberwithin{equation}{section}

\renewcommand{\bm}[1]{\vc{#1}}

\DeclareMathOperator*{\argmin}{arg\,min}
\def\avg#1{\left\{\mskip-5mu\left\{#1\right\}\mskip-5mu\right\}}

\def\div{\operatorname{div}}
\def\dt{\prtl_t}

\def\eps{\varepsilon}
\newcommand{\fdiv}{\overline{\div}}

\def\ff{\vc f}
\newcommand{\fgrad}{\overline{\nabla}} 
\newcommand{\feps}{\overline\eps}
\def\Hdiv{\vc H_{\div}}

\def\nn{\vc n}
\def\nnu{\boldsymbol\nu}

\def\prtl{\partial}

\def\Real{{\mathbb R}} 
\def\sigmapor{{\vc{\widehat\sigma}}}
\def\tn#1{{\mathbf{#1}}}    

\def\uu{\vc u}

\def\vc#1{\mathbf{\boldsymbol{#1}}} 

\def\vv{\vc v}

\def\ww{{\vc w}}
\def\xx{\vc x}

\def\zz{{\vc z}}

\let\oldequation\equation
\let\oldendequation\endequation

\renewenvironment{equation}{\begin{linenomath}\oldequation}{\oldendequation\end{linenomath}\par\nobreak\noindent}
\catcode`\*=11
\let\oldequationstar\equation*
\let\oldendequationstar\endequation*
\catcode`\*=12
\renewenvironment{equation*}{\begin{linenomath}\oldequationstar}{\oldendequationstar\end{linenomath}\par\nobreak\noindent}

\newcommand{\eq}[1]{\begin{equation}#1\end{equation}}
\newcommand{\eqs}[1]{\begin{equation*}#1\end{equation*}}
\newcommand{\ml}[1]{\begin{multline}#1\end{multline}}
\newcommand{\mls}[1]{\begin{multline*}#1\end{multline*}}

\begin{document}

\begin{frontmatter}



\title{On the parallel solution of hydro-mechanical problems with fracture networks and contact conditions}


\author[label1]{Jan Stebel\corref{cor1}}
\ead{jan.stebel@tul.cz}
\author[label2,label3]{Jakub Kružík}
\author[label2,label3]{David Horák}
\author[label1]{Jan Březina}
\author[label2,label3]{Michal Béreš}
\cortext[cor1]{Corresponding author.}

\affiliation[label1]{organization={Institute of New Technologies and Applied Informatics, Faculty of Mechatronics, Informatics and Interdisciplinary Studies, Technical University of~Liberec},
            addressline={Studentská~1402/2}, 
            city={Liberec},
            postcode={461 17}, 
            country={Czech Republic}}
\affiliation[label2]{organization={Institute of Geonics of the Czech Academy~of Sciences},
addressline={Studentská~1768}, 
city={Ostrava},
postcode={708~00}, 
country={Czech Republic}}
\affiliation[label3]{organization={Department of Applied Mathematics, VSB - Technical University of~Ostrava},
addressline={17.~listopadu~15/2172}, 
city={Ostrava},
postcode={708 00}, 
country={Czech Republic}}

\begin{abstract}
The paper presents a numerical method for simulating flow and mecha\-nics in fractured rock. The governing equations that couple the effects in the rock mass and in the fractures are obtained using the discrete fracture-matrix approach. The fracture flow is driven by the cubic law, and the contact conditions prevent fractures from self-penetration. A stable finite element discretization is proposed for the displacement-pressure-flux formulation. The resulting nonlinear algebraic system of equations and inequalities is decoupled using a robust iterative splitting into the linearized flow subproblem, and the quadratic programming problem for the mechanical part. The non-penetration conditions are solved by means of dualization and an optimal quadratic programming algorithm. The capability of the numerical scheme is demonstrated on a benchmark problem for \rev{tunnel} excavation with hundreds of fractures in 3D.
The paper's novelty consists in a combination of three crucial in\-gre\-dients: (i) application of discrete fracture-matrix approach to poroelasticity, (ii) robust iterative splitting of resulting nonlinear algebraic system working for real-world 3D problems, and (iii) efficient solution of its mechanical quadratic programming part with a large number of fractures in mutual contact by means of own solvers implemented into an in-house software library.
\end{abstract}




\end{frontmatter}
\nolinenumbers

\section{Introduction}

The interaction between fluid flow and mechanics in brittle rocks is an important phenomenon in many geotechnical applications. This type of rock is considered e.g. for designing deep geological repositories for radioactive waste\rev{, CO${}_2$ storage facilities or enhanced geothermal systems}. Numerical modelling of related problems is a complex task involving coupled physical processes at different geometrical scales, see \cite{Rutqvist2003rolea}. For the accurate resolution of phenomena in the rock matrix as well as in fractures and fault zones, the discrete fracture-matrix (DFM) approach has been introduced, see \cite{Berre2019Flowb}, as an alternative to traditional continuum-based and discrete fracture network (DFN) approaches. DFM models describe the physical processes both in the rock matrix and in the fracture network, as well as their interaction. The challenging task, however, is the discretization of geometries combining entities of different dimensions and the solution of large computational problems. Therefore, the number of fractures that have been used in DFM simulations so far is much smaller (typically units to tens) compared to DFN models (thousands), see e.g. \cite{berrone2015parallel,de2013synthetic,bergamaschi2023}.

\rev{
DFM models were originally developed for porous media flow in the pioneering works \cite{alboin2002modeling,flauraud2003,karimi2004,Martin2005Modeling}, and studied by many authors since then, see e.g. \cite{angot2009asymptotic,formaggia2014reduced,schwenck2015dimensionally,flemisch2016review}.
The coupling of DFM flow to rock mechanics was later considered e.g. in \cite{girault2015lubrication, girault2019mixed, bukac2017dimensional, hanowski2018hydromechanical}, omitting the mechanical response of fractures.
%
Models including fracture contact mechanics, which prevents self-penetration, were studied in \cite{garipov2016discrete,franceschini2020algebraically,Berge2020Finitea,blaheta2020bayesian,bonaldi2022numerical} and several other works, considering various types of friction laws or frictionless contact.
Solution of contact problems is a difficult task due to their nonlinearity and non-smoothness. Stable discretization of fracture contact mechanics was proposed in \cite{wohlmuth2011variationally}; see also \cite{franceschini2016novel} for the extension to plasticity.
We also mention modelling of hydraulic fracturing, where either the DFM models \cite{settgast2017fully} or alternative approaches have been used, such as the XFEM \cite{khoei2018enriched}, GFEM \cite{meschke2015generalized,shauer2020generalized} or the levelset method \cite{wheeler2014augmented}.
}

\rev{
Fractures in natural rocks are never smooth and void. In particular, due to healing, the fracture space can be filled by a solid material \cite{vass}. This leads us to propose a new DFM model combining the poroelastic rock matrix with poroelastic fractures. The DFM model extends the linear case from our previous work \cite{brst_mixed_hm} and is obtained by a semi-discretization process from physical principles, starting from the Biot equations \cite{Biot1941Generala}. The deformation of fractures is restricted by the minimal-width constraint, which acts as a contact condition, and the fracture permeability obeys the cubic law \cite{Snow1969Anisotropic}. The model is restricted to frictionless contact, and the fractures are assumed to be preexisting, i.e. we do not consider their initiation or propagation.
}

\rev{
Our second aim is to design a numerical approach for solving the coupled hydro-mechanical problem with a large number of fractures. That is important for simulations involving randomly generated fracture networks, for example in a Bayesian inversion analysis or stochastic homogenization. Similarly as in \cite{frigo2021efficient}, where the Biot model was solved in a domain without fractures, we use a finite element discretization based on the $P_1$ approximation of the displacement, \rev{the lowest-order Raviart-Thomas elements for the} velocity and \rev{$P_0$ approximation of the} pressure. We use compatible simplicial meshes for the rock and the fractures, which require a high-quality mesh generation tool. The model implementation however remains less complex in comparison with non-conforming or non-matching grids \cite{fumagalli2019conforming} or XFEM/GFEM type approaches \cite{fries2010extended}, that are advantageous e.g. for modelling fracture propagation. The resulting nonlinear algebraic system is decoupled using the fixed-stress splitting method \cite{Mikelic2013Convergence} to a sequence of steps involving the solution of a quadratic programming problem with inequality constraints. For their robust and efficient solution, we use the MPGP (Modified Proportioning with Gradient Projections) method \cite{Dos-book-09} implemented in the PERMON (Parallel, Efficient, Robust, Modular, Object-oriented, Numerical toolbox) library \cite{permonweb}. The numerical scalability of the method is demonstrated on a benchmark problem with hundreds of fractures.
}

The paper is organized as follows:
In \cref{sec:discrete_fract}, we present the hydro-mechanical model of fractured porous media. Then in \cref{sec:biot_discretization}, we describe the finite element approximation and the iterative decoupling of the resulting nonlinear problem to separate flow and mechanical subproblems. In \cref{sec:QP}, the solution of the quadratic programming problems is described.
\Cref{sec:numerical_examples} is devoted to \rev{a validation test and} a computational benchmark \rev{with} numerical experiments.

\section{Discrete fracture-matrix model of poroelasticity}\label{sec:discrete_fract}

This section describes the discrete fracture-matrix framework for fractured porous medium.
First, the semi-discrete differential operators are defined in \cref{sec:semidisc}. The hydro-mechanical model of the rock matrix based on the Biot theory of linear poroelasticity is introduced in \cref{sec:model_biot}. In \cref{sec:biot_dfm}, the semi-discrete operators are used to establish the equations of poroelasticity in fractures. Finally, the initial and boundary conditions are specified in \cref{sec:init_bc}. 

\subsection{Domain with fractures and semi-discrete operators}\label{sec:semidisc}

\begin{figure}
    \centering
    \includegraphics{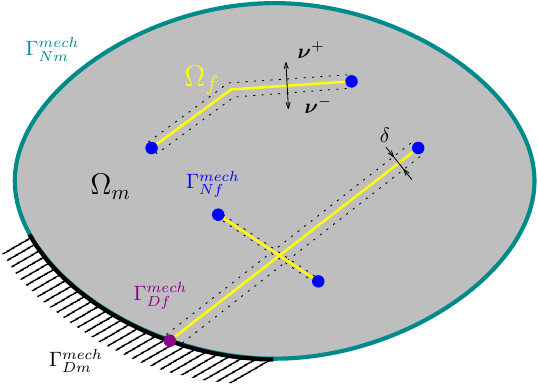}
    \caption{Illustration of a domain $\Omega_m$ representing the rock matrix with discrete fractures $\Omega_f$ and the boundaries, where Dirichlet and Neumann conditions are prescribed for the mechanics. Circles denote the boundary of $\Omega_f$. The flow boundary conditions can be prescribed on different parts of the boundary.}
    \label{fig:domain_dfn}
\end{figure}

Let $\Omega$ be a bounded domain in $\Real^d$, $d\in\{2,3\}$, with a Lipschitz boundary.
In this domain, we consider a set $\Omega_f\subset\Omega$, a union of $(d-1)$-dimensional polygons.
The set $\Omega_f$ represents the discrete fracture network, and the symbol $\Omega_m:=\Omega\setminus\overline\Omega_f$ defines the rock matrix; see \cref{fig:domain_dfn}.
With each discrete fracture, we associate a scalar field $\delta(\vc x)$, representing its initial \rev{width}.
If a point $\vc x\in\Omega_f$ belongs to a single fracture, we denote by $\nnu^+(\vc x)$ and $\nnu^-(\vc x):=-\nnu^+(\vc x)$ the unit normal vectors to $\Omega_f$ at $\vc x$.
Their orientation can be arbitrary but fixed for each fracture.

For a particular quantity, e.g. the pressure head $p$, we write $p=(p_m,p_f)$, where $p_m$ and $p_f$ is the value in $\Omega_m$ and $\Omega_f$, respectively.
The quantity $p$ can be discontinuous in the normal direction to the fracture; thus, we denote the values of $p_m$  on the opposite sides of the fracture by $p_m^+$ and $p_m^-$, respectively.

To simplify the notation of the equations in the fractures, we introduce the following operators (see \cite{brst_mixed_hm} for more details).
Let $\avg{f}:=\frac12(f^++f^-)$ denote the averaging operator on $\Omega_f$.
The semi-discrete gradient of a scalar field $p$ in $\Omega_f$ is
\eq{ \fgrad p := \nabla_\tau p_f + \avg{\Delta p\nnu}, }
where $\nabla_\tau$ stands for the tangential gradient in $\Omega_f$, and the normal gradient
is approximated by the average of differences across each side of the fracture:
\eq{ \rev{\avg{\Delta p\nnu} = \frac12(\Delta^+ p\nnu^+ + \Delta^- p\nnu^-),\quad} \Delta^\pm p := \frac2\delta(p_m^\pm-p_f). }
Next, we define the semi-discrete gradient on the positive and negative sides of the fracture:
\eq{ \fgrad^\pm p := \nabla_\tau p_f + \Delta^\pm p\nnu^\pm. }
%
In the same way, we introduce the semi-discrete gradient for a vector field $\uu=(\uu_m,\uu_f)$\rev{, where $\uu_m:\Omega_m\to\Real^d$ and $\uu_f:\Omega_f\to\Real^d$,} and its symmetric part:
\eq{ \fgrad\uu := \nabla_\tau\uu_f + \avg{\Delta\uu\otimes\nnu},\quad \feps(\uu):=\frac12(\fgrad\uu+(\fgrad\uu)^\top). }
Analogously, we can define the semi-discrete divergence of a vector and tensor quantity:
\eq{ \fdiv\uu := \div_\tau\uu_f + \avg{\Delta\uu\cdot\nnu}, \qquad \fdiv\vc\sigma := \div_\tau\vc\sigma_f + \avg{\Delta\vc\sigma \nnu}, }
and other operators, e.g. $\feps^\pm(\uu)$, $\fdiv^\pm\uu$ etc.
\rev{The tangential part of a vector $\uu$ will be denoted by $\uu_\tau:=\uu-(\uu\cdot\nnu^+)\nnu^+$.}

\subsection{Biot's poroelasticity model for rock matrix}\label{sec:model_biot}

The mechanical behaviour of the rock matrix is based on the equilibrium equation
\eq{\label{eq:el} -\div\rev{\sigmapor_m} = \vc f_m \mbox{ in }I\times\Omega_m,}%
where $\rev{\sigmapor_m:=\vc\sigma_m - \alpha_m\varrho g p_m\tn I}$ is the \rev{total} stress tensor, $\alpha_m$ is the Biot-Willis coefficient, $\varrho$ is the density of the liquid, $g$ is the gravitational acceleration, $p_m$ is the pressure head, $\vc f_m$ is the load and $I:=(0,T)$ is a time interval.
The mass balance in the rock is given by the equation
\eq{\label{eq:fl} \prtl_t(S_m p_m+\alpha_m\div\uu_m) + \div\vc v_m = s_m \mbox{ in }I\times \Omega_m, }%
where $S_m$ is the storativity, $\uu_m$ is the displacement, $\vc v_m$ is the Darcy flow velocity, and $s_m$ is the volume source of the fluid.
We consider Hooke's and Darcy's laws as the constitutive relations for the \rev{effective} stress tensor $\vc\sigma_m$ and velocity, respectively:
\eq{ \vc\sigma_m - \vc\sigma_{0m}=\tn C_m\eps(\uu_m), \qquad \vc v_m = -\tn K_m(\nabla p_m+\vc g). }
Here $\vc\sigma_{0m}$ is the initial stress, and $\vc g$ is the \rev{gradient of the vertical distance (typically $\vc g=\nabla x_d$)}.
For isotropic media, the Hooke tensor and the hydraulic conductivity tensor have the form:
\eq{\label{eq:isotropic_laws_m} \tn C_m = \frac{E_m}{1+\nu_m}\tn I_4 + \frac{E_m\nu_m}{(1+\nu_m)(1-2\nu_m)}\tn I\otimes\tn I,\qquad \tn K_m=k_m\tn I, }
where $E_m,\nu_m$ is the Young modulus and the Poisson ratio, respectively, $k_m$ is the (scalar) hydraulic conductivity, and $\tn I,\tn I_4$ are the second and fourth order unit tensors, respectively.

\subsection{Poroelastic model of fractures}\label{sec:biot_dfm}

\begin{figure}
    \centering
    \includegraphics[width=0.75\textwidth]{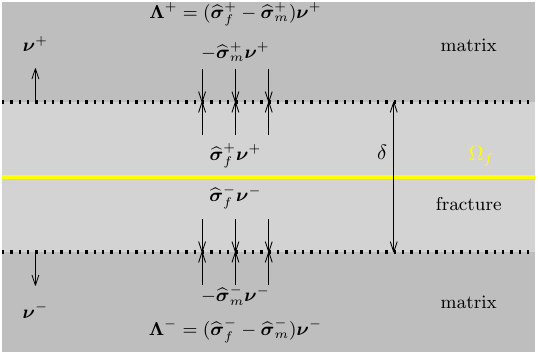}
    \caption{\rev{Scheme of contact forces acting on the matrix-fracture interface.}}
    \label{fig:contact-forces}
\end{figure}

The pressure head $p_f$ and the displacement $\uu_f$ in $\Omega_f$ are given by the semi-discretized equations
\begin{linenomath}\begin{align}
\label{eq:mixed_dim_problem_el_frac} -\fdiv\rev{(\delta\sigmapor)} &= \rev{\delta}\vc f_f &&\mbox{ in }I\times\Omega_f,\\
\label{eq:mixed_dim_problem_fl_frac} \prtl_t\left(\rev{\delta}S_f p_f + \rev{\delta}\alpha_f\fdiv\,\uu\right) +\fdiv\,\rev{(\delta\vv_m,\vv_f)} &= \rev{\delta}s_f &&\mbox{ in }I\times\Omega_f,
\end{align}\end{linenomath}
where, \rev{in accordance with the notation from \cref{sec:semidisc}, $\fdiv\,(\delta\vv_m,\vv_f)=\div_\tau\vv_f+\vv_m^+\cdot\nnu^+ + \vv_m^-\cdot\nnu^-$.} \rev{The total} stress tensor and the flux in the~fracture are defined by the relations
\begin{linenomath}\eq{
\label{eq:const_laws_f} 
\rev{\sigmapor_f = \vc\sigma_{0f}+\tn C_f\feps(\uu)-\alpha_f\varrho g p_f\tn I},
\qquad \vc v_f = -\rev{\delta\tn K_f(\nabla_\tau p_f+\vc g_\tau)}.
}\end{linenomath}
Here $\vc\sigma_{0f}$, $\tn C_f$ \rev{and $\alpha_f$} is the initial stress, the elasticity tensor \rev{and the Biot-Willis coefficient} in the fractures, respectively. \rev{The elastic tensor $\tn C_f$ is given by a similar relation as in \cref{eq:isotropic_laws_m} via coefficients $E_f$, $\nu_f$, which can be determined either from experimental data or rock joint constitutive models \cite{bandis1983fundamentals}.} The fracture permeability tensor $\tn K_f$ \rev{is} given by the cubic law \cite{Snow1969Anisotropic}:
\begin{linenomath}\begin{equation*}
     \tn K_f = k_f(\uu)\tn I,\qquad k_f(\uu) = \frac{\eta\varrho g}{12\mu} a_f(\uu)^2
\end{equation*}\end{linenomath}
in terms of the aperture $a_f = a_f(\uu) := \delta + \rev{\uu_m^+\cdot\nnu^++\uu_m^-\cdot\nnu^-}$, \rev{coefficient of} fracture roughness $\eta\in(0,1]$ and fluid viscosity $\mu$.
\rev{Other symbols i.e. $\ff_f$, $S_f$ and $s_f$ have the same physical meaning as their counterparts in $\Omega_m$.
To write down the matrix-fracture interface conditions, we introduce the normal stress and contact force on each fracture side:
\begin{linenomath}\begin{align}\label{eq:flux_stress_pm}
\vc \sigmapor_f^\pm&:=\vc\sigma_{0f}+\tn C_f\feps^\pm(\uu)-\alpha_f\varrho g p_f\tn I,\\
\vc\Lambda^\pm&:=(\sigmapor_f^\pm-\sigmapor_m^\pm)\nnu^\pm.
\end{align}\end{linenomath}
Here, $\vc\Lambda^\pm$ is the force acting on the fracture side as a result of the tractions from the matrix and the fracture, see \cref{fig:contact-forces}.
On $\Omega_f$ we impose:}
\begin{itemize}
    \item Continuity of the flux:
    \eq{\label{eq:cont_flux} \vc v_m^\pm\cdot\nnu^\pm = \rev{-\tn K_f(\fgrad^\pm p+\vc g)\cdot\nnu^\pm};}
    \item \rev{Balance of contact forces on fracture sides:}
    \eq{\label{eq:contact_balance}\rev{\vc\Lambda^+= -\vc\Lambda^-;}}
    \item \rev{Zero friction:}
    \eq{\label{eq:contact_friction}\rev{\vc\Lambda_\tau^\pm = \vc 0;}}
    \item \rev{Non-penetration:
    \eq{ a_f\ge \delta_{min}}
    with $\delta_{min}>0$ the minimal \rev{width} of the fracture;}
    \item \rev{Complementarity relations:
    \begin{linenomath}\begin{align}    
    \label{eq:contact_stress}\vc\Lambda^\pm\cdot\nnu^\pm &\ge 0,\\
    \label{eq:contact_complementarity}\avg{\vc\Lambda\cdot\nnu}(a_f-\delta_{min}) &= 0.
    \end{align}\end{linenomath}
    These relations imply that the fracture closes to its minimal width when contact forces are nonzero and, vice versa, contact forces vanish when the fracture is open ($a_f>\delta_{min}$).}
\end{itemize}
At intersections of fractures, we require the displacement $\uu_f$ and pressure $p_f$ to be continuous.

\subsection{Initial and boundary conditions}\label{sec:init_bc}

The DFM Biot model is supplemented with the initial condition for the pressure head
\eq{\label{eq:init_bc} p_{m|t=0} = p_{0m} \mbox{ in }\Omega_m,\qquad p_{f|t=0} = p_{0f}\mbox{ in }\Omega_f. }
The initial displacement is determined by the initial pressure head through the elasticity \cref{eq:el,eq:mixed_dim_problem_el_frac}, respectively.
For the flow equations, we consider the boundary conditions 
\begin{linenomath}\begin{align}
\label{eq:bc_flow_D}p_* &= p_{D*} &&\mbox{ on } I\times\Gamma_{D*}^{flow},\\
\label{eq:bc_flow_N}\vc v_*\cdot\vc n &= v_{N*} &&\mbox{ on }I\times\Gamma_{N*}^{flow},~*\in\{m,f\},
\end{align}\end{linenomath}
where $p_D:=(p_{Dm},p_{Df})$ is a given pressure head, $v_N:=(v_{Nm},v_{Nf})$ a given flux, $\nn$ is the unit outward normal vector, and $\Gamma_{D*}^{flow}$, $\Gamma_{N*}^{flow}$ denote the corresponding parts of the boundary:
\[ \partial\Omega_* = \overline\Gamma_{D*}^{flow} \cup \overline\Gamma_{N*}^{flow},~*\in\{m,f\}. \]
Further, we consider the Dirichlet and Neumann boundary conditions for the mechanics:
\begin{linenomath}\begin{align}
\label{eq:bc_mech_D} \uu_* &= \uu_{D*} &&\mbox{ on } I\times\Gamma_{D*}^{mech},\\
\label{eq:bc_mech_N}\rev{\sigmapor_*}\vc n &= \vc t_{N*} &&\mbox{ on }I\times\Gamma_{N*}^{mech},~*\in\{m,f\}
\end{align}\end{linenomath}
with given displacement $\uu_D:=(\uu_{Dm},\uu_{Df})$ and surface traction $\vc t_N:=(\vc t_{Nm},\vc t_{Nf})$.
Here again $\Gamma_{D*}^{mech}$, $\Gamma_{N*}^{mech}$ are the complementary parts of the boundary such that
\[ \partial\Omega_* = \overline\Gamma_{D*}^{mech}\cup\overline\Gamma_{N*}^{mech},~*\in\{m,f\}. \]

\rev{Altogether, \crefrange{eq:el}{eq:bc_mech_N} define the coupled DFM hydro-mechanical problem that will be referred to as Problem (DFM-HM).}

\section{Finite element approximation and algebraic formulation} \label{sec:biot_discretization}

\begin{figure}
    \centering
    \begin{subfigure}[b]{0.49\textwidth}
    \includegraphics[width=\textwidth]{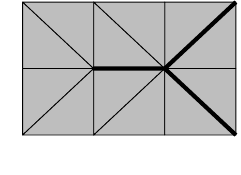}
    \caption{Compatible mesh of rock matrix and discrete fractures.}
    \end{subfigure}
    \begin{subfigure}[b]{0.49\textwidth}
    \includegraphics[width=\textwidth]{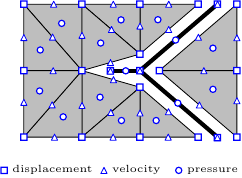}
    \caption{Distribution of degrees of freedom for the $P_1/RT_0$ finite elements.}
    \end{subfigure}
    \caption{Scheme of spacial discretization in 2D.}
    \label{fig:mesh-dofs}
\end{figure}

\rev{Problem (DFM-HM)} is discretized by the finite element method and the implicit Euler time-stepping.
For discretizing the rock matrix $\Omega_m$ and the fracture network $\Omega_f$, we use simplicial meshes denoted $\mathcal T_m^h$ and $\mathcal T_f^h$, respectively.
We require compatibility of these meshes, which means that the elements of $\mathcal T_f^h$ coincide with the faces of adjacent elements of $\mathcal T_m^h$.
On top of $(\mathcal T_m^h,\mathcal T_f^h)$ we define the following finite element spaces:
\begin{linenomath}
\mls{ \mathcal Z^h:=\{(\zz_m^h,\zz_f^h)\in H^1(\Omega_m;\Real^d)\times H^1(\Omega_f;\Real^d);\\
    \zz_{*|E}^h\in P_1(E;\Real^d)~\forall E\in\mathcal T_*^h,~*\in\{m,f\}\}, } \end{linenomath}
\[ \mathcal Q^h:=\{(q_m^h,q_f^h)\in L^2(\Omega_m)\times L^2(\Omega_f);q_{*|E}^h\in P_0(E)~\forall E\in\mathcal T_*^h~*\in\{m,f\}\}, \]
\begin{linenomath}\mls{ \mathcal W^h(v_N):=\{(\ww_m^h,\ww_f^h)\in \Hdiv(\Omega_m)\times\Hdiv(\Omega_f);\\
     \ww_{*|E}^h\in RT_0(E)~\forall E\in\mathcal T_*^h,~
     \ww_*^h\cdot\nn=\Pi^h_0 (v_{N*})\mbox{ on }\Gamma_{N*}^{flow},~*\in\{m,f\} \}, }\end{linenomath}
where $P_k$, $RT_k$, and $\Pi^h_k$ denote the space of polynomials of degree up to $k$, the $k$-th order Raviart-Thomas space, and the projector onto piecewise polynomials of degree up to $k$, \rev{respectively}.
We also introduce the cone of admissible displacements:
\begin{linenomath}\mls{
    \mathcal K^h := \{\zz^h\in\mathcal Z^h;~\zz^h_*=\Pi^h_1(\uu_{D*}) \mbox{ on }\Gamma_{D*}^{mech},~*\in\{m,f\},\\
    ~\forall E\in\mathcal T_f^h:~\tfrac1{|E|}\int_E a_f(\zz^h)\ge\delta_{min} \}.
}\end{linenomath}
We refer to \cref{fig:mesh-dofs} for the illustration of the computational mesh and distribution of degrees of freedom for finite element functions.
\rev{We point out that the discrete functions from $\mathcal Z^h$, $\mathcal Q^h$, and $\mathcal W^h(v_N)$ are generally discontinuous on the matrix-fracture interface but continuous on fracture intersections. The continuity of the flux is enforced by the discretized equations given below.}

\subsection{Discretized problem}
 \rev{The discretized version of the Problem (DFM-HM) stems from the weak formulation (we refer to \cref{sec:appendix_weak} for its derivation).}
In each discrete time $t_n:=n\Delta t$, we solve the following nonlinear problem:
Find $(\uu^{h,n},p^{h,n},\vv^{h,n})\in\mathcal K^h\times\mathcal Q^h\times\mathcal W^h(v_N)$ such that for all $(\zz^h,q^h,\ww^h)\in\mathcal K^h\times\mathcal Q^h\times\mathcal W^h(0)$:
\begin{linenomath}\begin{align*}
a(\uu^{h,n},\zz^h-\uu^{h,n}) - b(\zz^h-\uu^{h,n},p^{h,n}) &\rev{\ge} f_1(\zz^h-\uu^{h,n}),\\
c(p^{h,n},q^h) + b(\uu^{h,n},q^h) + \Delta t d(\vv^{h,n},q^h) &= \Delta t f_2(q^h) + c(p^{h,n-1},q^h) + b(\uu^{h,n-1},q^h),\\
e(\uu^{h,n};\vv^{h,n},\ww^h) - d(\ww^h,p^{h,n}) &= f_3(\ww^h),
\end{align*}\end{linenomath}
where the forms are defined in \cref{eq:weak_forms}.

\subsection{Iterative splitting}
The solution of the discretized problem is obtained using the fixed-stress splitting method, which leads to the following iterative process for each time $t_n$:
\begin{itemize}
    \item Set $p^{h,n,0}:=p^{h,n-1}$;
    \item For $i=1,2,...,i_{max}$
    \begin{itemize}
        \item[(i)] Find $\uu^{h,n,i}\in\mathcal K^h$ such that for all $\zz^h\in\mathcal K^h$:
        \[ a(\uu^{h,n,i},\zz^h-\uu^{h,n,i}) - b(\zz^h-\uu^{h,n,i},p^{h,n,i-1}) \rev{\ge} f_1(\zz^h-\uu^{h,n,i}); \]
        \item[(ii)] Find $(p^{h,n,i},\vv^{h,n,i})\in\mathcal Q^h\times\mathcal W^h(v_N)$ such that for all $(q^h,\ww^h)\in\mathcal Q^h\times \mathcal W^h$:
        \mls{ c(p^{h,n,i},q^h) + c_\beta(p^{h,n,i},q^h) + b(\uu^{h,n,i},q^h) + \Delta t d(\vv^{h,n,i},q^h)\\
        = \Delta t f_2(q^h) + c(p^{h,n-1},q^h) + c_\beta(p^{h,n,i-1},q^h) + b(\uu^{h,n-1},q^h), }
        \[ e(\uu^{h,n,i};\vv^{h,n,i},\ww^h) - d(\ww^h,p^{h,n,i}) = f_3(\ww^h). \]
        \item[(iii)] \rev{If $\|p^{h,n,i}-p^{h,n,i-1}\|\le Tol\|p^{h,n,i-1}\|$ then set $(\uu^{h,n},p^{h,n},\vv^{h,n}):=(\uu^{h,n,i},p^{h,n,i},\vv^{h,n,i})$ and stop.}
    \end{itemize}
\end{itemize}
The bilinear form
\[ c_\beta(p,q):=\int_{\Omega_m}\beta_m p_m q_m + \int_{\Omega_f}\rev{\delta}\beta_f p_f q_f, \]
used in step (ii), represents a stabilization term, \rev{$i_{max}$ is the maximal number of iterations and $Tol$ is the relative tolerance for the stopping criterion}.
In accordance with theoretical studies (\cite{brst_mixed_hm}, see also \cite{both2017robust}), we set
\[ \beta_* := \frac{\alpha_*^2(1+\nu_*)(1-2\nu_*)}{E_*},~*\in\{m,f\}, \]
which ensures convergence of the iterative scheme for the linear Biot system.

\subsection{Algebraic representation}
Steps (i)-(ii) of the iterative scheme can be written in the algebraic form as follows:
Find the unknown vectors of degrees of freedom $(\uu_{n,i},\vc p_{n,i},\vv_{n,i})$ such that
\begin{subequations}\label{eq:hm_alg}
\eq{ \label{eq:mech_alg_split} \argmin_{\uu_{n,i}}\left(\frac12\tn A\uu_{n,i}-\vc f_{1,n,i}\right) \cdot \uu_{n,i} \mbox{ subject to }\tn B_I\uu_{n,i}\le\vc c_I, }
\eq{ \label{eq:flow_alg_split} \begin{bmatrix}\tn C+\tn C_\beta & \Delta t\tn D^\top\\-\tn D & \tn E(\uu_{n,i})\end{bmatrix} \begin{bmatrix}\vc p_{n,i}\\\vv_{n,i}\end{bmatrix} = \begin{bmatrix}\vc f_{2,n,i}\\\vc f_3\end{bmatrix}, }
where the matrices $\tn A,...,\tn E$ and vectors $\vc f_1,...,\vc f_3$ 
arise from the corresponding forms $a,...,e,f_1,...,f_3$.
Further
\[ \vc f_{1,n,i}:=\vc f_1+\tn B^\top\vc p_{n,i-1}, \quad \vc f_{2,n,i}:=\Delta t\vc f_2 + \tn C\vc p_{n-1} + \tn C_\beta\vc p_{n,i-1} - \tn B(\uu_{n,i}-\uu_{n-1}), \]
and $\tn B_I$, $\vc c_I$ is the matrix and vector of constraints realizing the contact conditions.
\end{subequations}

We observe that within the above iterative process, the hydraulic subproblem reduces to the linear algebraic system \cref{eq:flow_alg_split}. This system is further reduced using the hybridization technique to a smaller problem for the pressure traces on element faces \rev{(see \cite{sistek2015bddc} for more details)}.
The matrix of the latter problem is symmetric positive definite, i.e. the problem can be solved e.g. by the conjugate gradient method.
The mechanical problem \cref{eq:mech_alg_split} is a quadratic programming (QP) problem with linear inequality constraints, whose accurate and efficient solution is a key part of the numerical method.
In the following section, we describe its solution.

\section{Solution of mechanical QP problems}\label{sec:QP}

From the mathematical point of view, the original mechanical primal problem \cref{eq:mech_alg_split} is a QP problem of the following type
\begin{equation}
 \argmin_{{\vc x}}\frac{1}{2}\vc x^T\tn A \vc x-\vc b^T\vc x\;\;\text{s.t.}\;\;\tn B_I\vc x\leq \vc c_I,\label{QP_prim}
\end{equation}
where $\tn A$ is a symmetric positive definite matrix, $\vc b := \vc f_{1,n,i}$ and $\vc x := \uu_{n,i}$.

In our case, the original QP problem \cref{QP_prim} is transformed into the dual one, obtaining a QP problem with a simple bound 
\begin{equation}
 \argmin_{{\vc {\lambda}}}\frac{1}{2}\vc {\lambda}^T\tn F \vc {\lambda}-\vc d^T\vc {\lambda}\;\;\text{s.t.}\;\;\vc {\lambda}\geq \rev{\vc 0} ,\label{QP_simplebound}
\end{equation}
where $\tn{F}:=\tn B_I\tn{A}^{-1}\tn B_I^{T}$ and $\vc d := \tn B_I\bm{A}^{-1}\vc{b} - \vc c_I$.

QP problem \cref{QP_simplebound} can be efficiently solved by the MPGP algorithm \cite{Dos-book-09}.

Once $\vc{\lambda}$ is found, the solution $\vc{x}$ can be evaluated
by the formula
\begin{equation}
\vc{x}=\tn{A}^{-1}(\vc{b}-\tn B_I^{T}\vc{\lambda}).\label{eq3}
\end{equation}
More details can be found in \cite{HapHorPos-HPCSE-16}.

\section{Numerical experiments}\label{sec:numerical_examples}

In this section, we describe the implementation of the numerical scheme presented above.
We introduce the software codes used for the subsequent numerical simulations on computer clusters.
\rev{Then we present a 2D validation test.}
Finally, we show the solution of a 3D computational benchmark and perform numerical experiments with a parallel solution.

\subsection{Software implementation}

The coupled hydro-mechanical model has been implemented in the open-source code Flow123d \cite{Brezina2021Flow123d}, \rev{which is} \rev{developed by the authors}.
The software provides the solution of finite element models of processes in porous media built on top of meshes, that combine tetrahedral, triangular, and line elements to represent the bulk material, fractures, and their intersections, respectively. The implementation supports running on parallel architectures. The computational mesh is \rev{created by the GMSH mesh generator \cite{geuzaine2009gmsh} and} partitioned into subdomains using the ParMETIS \cite{karypis1998parallel} library. Parallel matrices and vectors are assembled in the PETSc format \cite{petsc-www,petsc-user-ref}.
We use PETSc solvers to solve linear algebra problems. Specifically, for hydraulic problems, we utilize the parallel version of conjugate gradients combined with the BoomerAMG algebraic multigrid preconditioner  \cite{yang2002boomeramg}.

For the parallel solution of \rev{the mechanical QP problems}, we used the software package PERMON \cite{permonweb,HapHorPos-HPCSE-16}, \rev{which is} \rev{developed by the authors}. 
PERMON is built on top of PETSc, mainly its linear algebraic part, and it is available for free under the FreeBSD open-source license \cite{permongit}.

PERMON extends PETSc by adding specific functionality and algorithms for large-scale sparse QP problems, as well as FETI-type domain decomposition methods \cite{DosHorKuc-CNME-06}. The same coding style is employed, allowing users familiar with PETSc to utilize these features with minimal effort. Among the primary applications are contact problems in mechanics, such as those involving opening and closing fractures, like the problem represented by  \cref{eq:mech_alg_split}.

The main \rev{part of PERMON,} PermonQP, provides a base for solving QP problems. It includes data structures, transformations, algorithms, and supporting functions for QP. 
The parallelization is primarily done through the row-wise distribution of the matrices and vectors with communication realized by an MPI library.

The action of the matrix inverse $\bm{A}^{-1}$ \rev{needed in \cref{QP_simplebound,eq3}} was done by the Cholesky factorization of $\bm{A}$ implemented in the SuperLU\_DIST library \cite{SuperLUDist-AMC-03}




\rev{
\subsection{2D validation test}
To validate the coupled hydro-mechanical model with non-linear hydraulic conductivity and contact conditions on fractures, we choose a problem of fluid injection into a one-dimensional fracture located in the middle of a planar impermeable block.
This problem has been used to validate many computational codes such as ROCMAS II \cite{noorishad1992theoretical}, FEHM \cite{bower1997numerical}, and COMPASS \cite{chen2020numerical}. The setting of the test is chosen to mimic the semi-analytical similarity solution from \cite{Wijesinghe1986},
parameter values are taken from \cite{watanabe2012lower}.
The computational domain is the rectangle $(0, 25) \times (0, 1)$, with the horizontal fracture $(0, 25) \times \{ \tfrac 12 \}$.
Initially, there is a uniform fluid pressure $p_0 = 11$ MPa in the whole domain.
A sufficiently large vertical initial stress of 50 MPa is applied in order to keep the fracture width at the minimal value $\delta_{min}=10^{-5}$ m.
During the test, the fluid is injected at the left fracture end under the pressure $p_{Df} = 11.9$ MPa, while at the right end, the pressure is kept at the initial value $p_0$.
The values of hydraulic and mechanical parameters are listed in \cref{tab:injection-params}.
\begin{table}
    \centering
    \caption{\rev{Physical parameters of the fluid injection test problem.}}
    \label{tab:injection-params}
    \begin{tabular}{l|c|r|r|c}
    \multicolumn{2}{c}{Quantity} & \multicolumn{1}{c}{Rock} & \multicolumn{1}{c}{Fracture}& Unit\\
    \hline
    hydraulic conductivity & $k_m$ & $9.81\cdot10^{-15}$ & $-$ & m/s\\
    storativity & $S_m,S_f$ & $9.81\cdot10^{-7}$ & $0$ & 1/m\\
    Young modulus & $E_m$ & $60$ & $\kappa\delta$ & GPa\\
    joint stiffness & $\kappa$ & $-$ & 100 & GPa/m\\
    Poisson's ratio & $\nu_m,\nu_f$ & 0 & 0 & $-$\\
    roughness coefficient & $\eta$ & $-$ & 1 & $-$\\
    minimal fracture width & $\delta_{min}$ & $-$ & $10^{-5}$ & m\\
    Biot-Willis coefficient & $\alpha_m,\alpha_f$ & $1$ & $1$ & $-$\\
    fluid density & $\varrho$ & $10^3$ & $10^3$ & kg/m${}^3$\\
    fluid viscosity & $\mu$ & $10^{-3}$ & $10^{-3}$ & Pa$\cdot$s\\
    \end{tabular}
\end{table}
In particular, the fracture hydraulic conductivity is given by the cubic law.
In this test, we consider the fracture width $\delta$ to be solution-dependent, i.e. $\delta=a_f=\delta_{min}+\uu_m^+\cdot\nnu^++\uu_m^-\cdot\nnu^-$.
The fracture Young modulus is set to $E_f=\kappa\delta$, where $\kappa=100$ GPa/m is the fracture joint stiffness.
The boundary of the rock domain is impermeable.
On the bottom side of the domain, we set zero displacement condition. The left and right sides are subject to the roller boundary condition $\uu_m\cdot\nn = 0$. 
The geometry and the boundary conditions are depicted in \cref{fig:injection_test_geom}.
\begin{figure}
    \includegraphics[width=\textwidth]{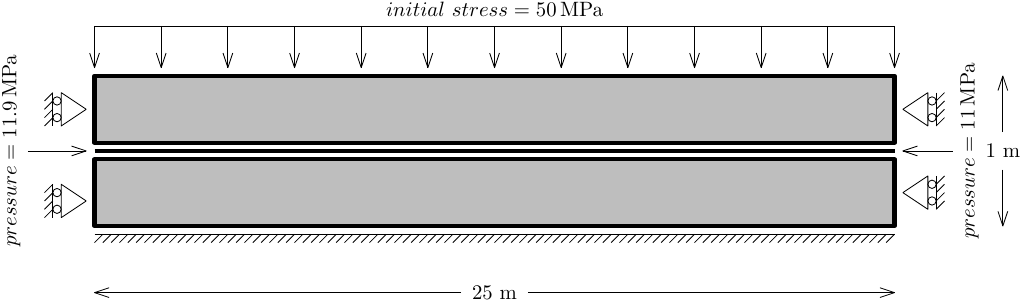}
    \caption{\rev{Geometry, initial and boundary conditions of the fluid injection test.}}
    \label{fig:injection_test_geom}
\end{figure}
}

\rev{
The problem was computed on a uniform mesh consisting of 22,000 triangles and 500 fracture line elements. 
The computational time interval is [0,~2000 s] with the time step $\Delta t = 1$ s. The reference semi-analytical solution from \cite{Wijesinghe1986}, to which the numerical results were compared, was obtained by solving a system of nonlinear ordinary differential equations with an implicitly given initial condition. We used the MATLAB function \texttt{ode45} and the shooting method for this task.
The progress of fracture opening at selected times is depicted in \cref{fig:injection_test_apertures}.
The numerical results exhibit an excellent agreement with the semi-analytical solution.
\begin{figure}
    \includegraphics[width=\textwidth]{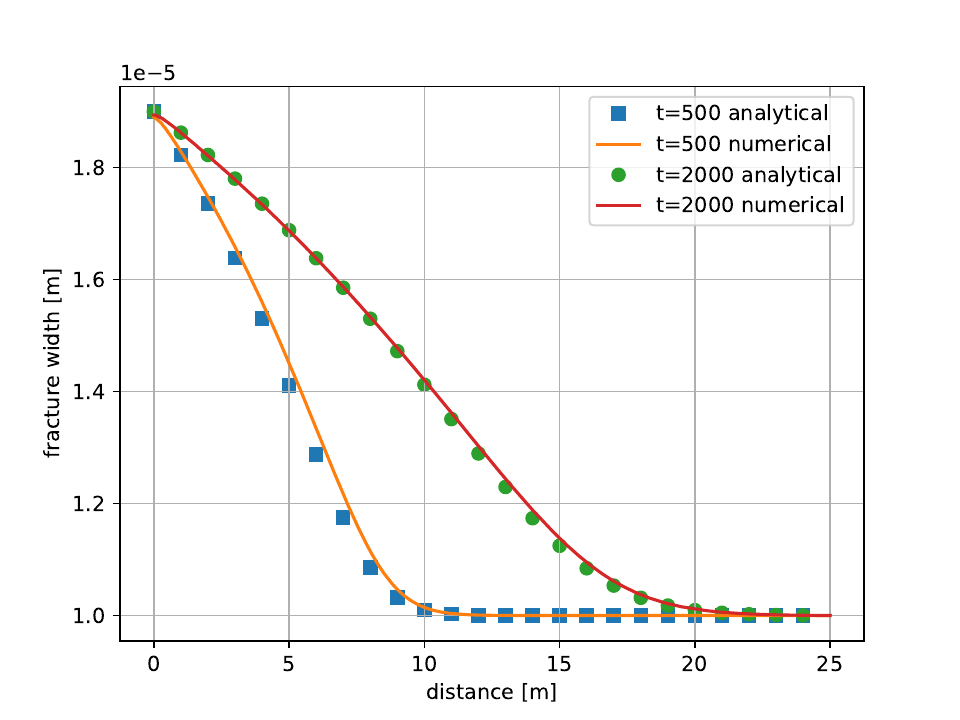}
    \caption{\rev{Fluid injection test: Comparison of fracture apertures obtained by numerical solution to the semi-analytical results at times 500 s and 2000 s.}}
    \label{fig:injection_test_apertures}
\end{figure}
}

\subsection{Benchmark - rock relaxation during \rev{tunnel} excavation}

\begin{figure}
    \centering
    \begin{subfigure}[b]{\textwidth}\centering
    \includegraphics[width=0.7\textwidth]{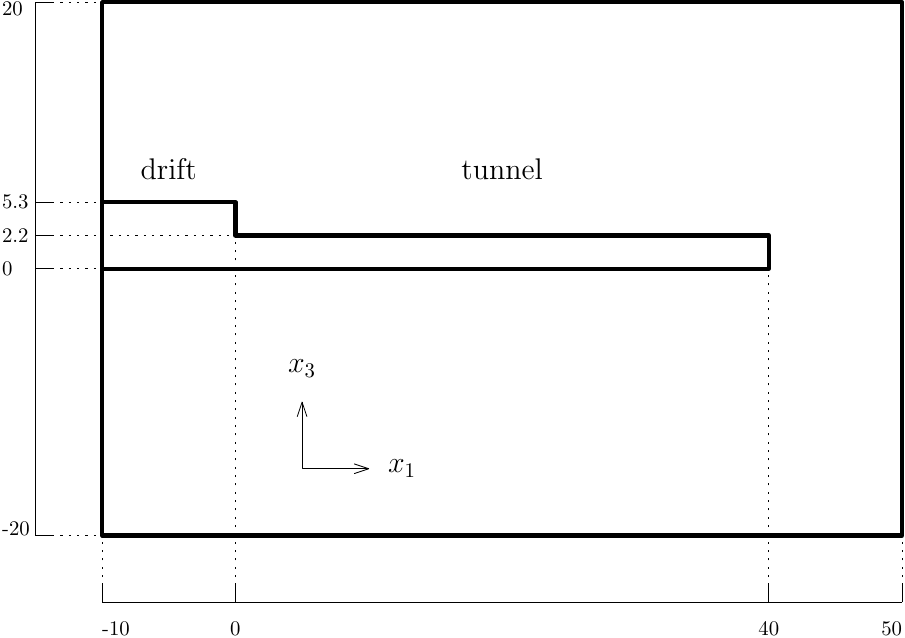}
    \caption{Side view of the computational domain.}
    \end{subfigure}
    \hspace{\stretch{1}}
    \begin{subfigure}[b]{\textwidth}
    \includegraphics[width=\textwidth]{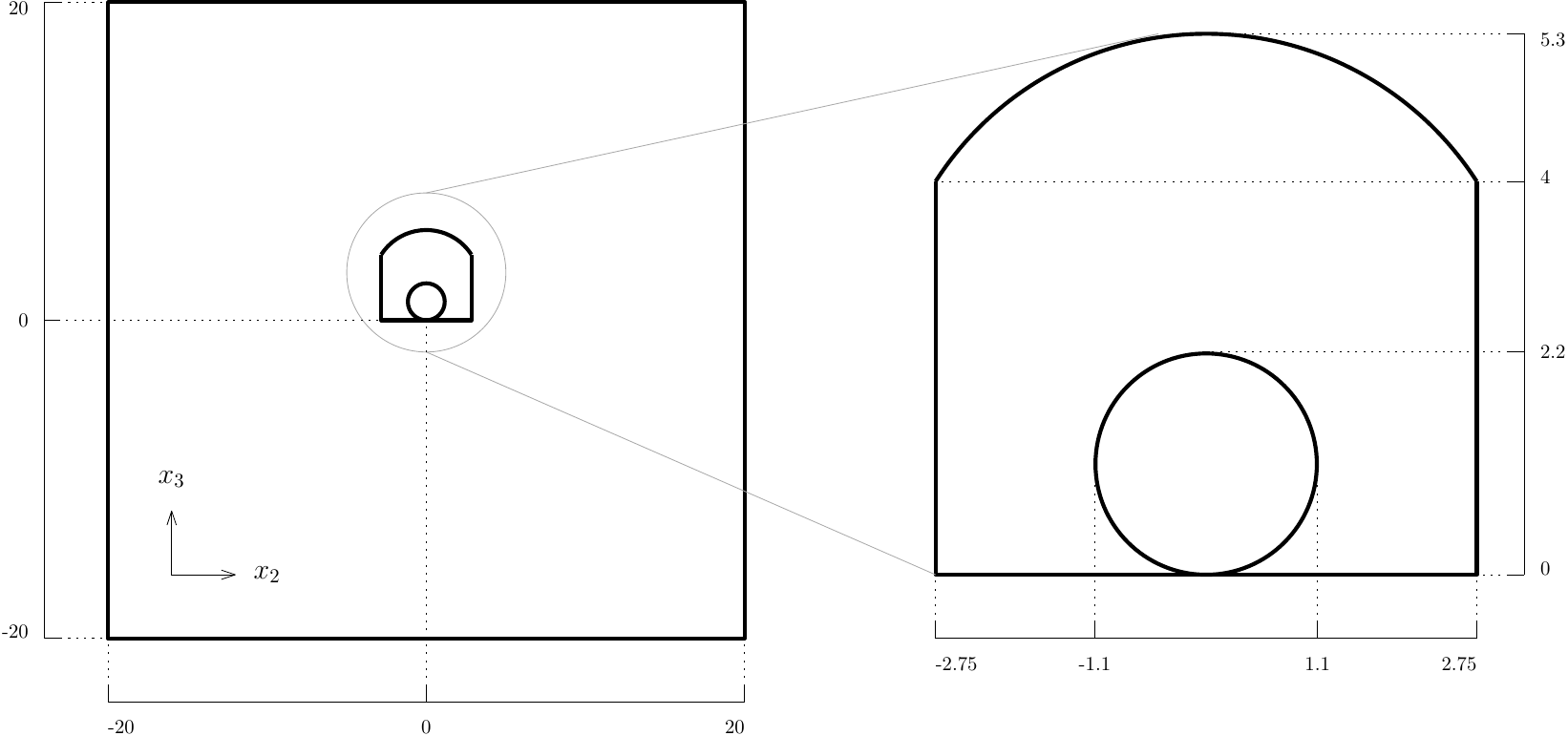}
    \caption{\rev{Front view of the computational domain; drift and tunnel dimensions}.}
    \end{subfigure}
    \caption{The geometry of the \rev{tunnel} excavation test problem.}
    \label{fig:tunnel-geo}
\end{figure}

We demonstrate the use of the presented numerical method in a realistic 3D geometry with a discrete fracture network.
The model simulates the excavation of a cylindrical \rev{tunnel} extending a \rev{large-profile drift} and, more importantly, the subsequent changes in rock pressure resulting from this excavation. It is assumed that the highly permeable fractures affect the pressure evolution in the domain. We compare the impact of two randomly generated discrete fracture networks with 200 and 400 fractures, respectively.
The computational domain $\Omega_m$ has a fixed geometry given by the block $(-10,50)\times(-20,20)\times(-20,20)$ m from which the \rev{drift} and the 40 m long cylindrical \rev{tunnel} \rev{are cut out}.
The dimensions of the domain are given in \cref{fig:tunnel-geo}.

In order to describe the boundary conditions, we define the following sets:
\begin{linenomath}\begin{align*}
    \Gamma^o &:= \{\xx\in\partial\Omega_m;~x_1=50 \mbox{ or }x_2=\pm20 \mbox{ or }x_3=\pm20\}, && \mbox{ (outer boundary)}\\
    \Gamma^n &:= \{\xx\in\partial\Omega_m;~x_1=-10\}, && \mbox{ (symmetry plane)}\\
    \Gamma^{\rev{d}} &:= \{\xx\in\partial\Omega_m\setminus\Gamma^o;~-10<x_1\le0\}, && \mbox{ (\rev{drift} surface)}\\
    \Gamma^{\rev{t}} &:= \{\xx\in\partial\Omega_m\setminus\Gamma^o;~0<x_1\le40\}. && \mbox{ (\rev{tunnel} surface)}
\end{align*}\end{linenomath}
We prescribe the following boundary conditions:
\begin{linenomath}\begin{align*}
    &\text{flow}&&\text{\hspace{-2em}mechanics}&&\\
    p_D&=300,& \uu\cdot\nn&=0 && \mbox{ on }\Gamma^o,\\
    p_D&=0,& \vc t_N&=\vc0 && \mbox{ on }\Gamma^{\rev{d}},\\
    \rev{v}_N&=0,& \uu\cdot\nn&=0 && \mbox{ on }\Gamma^n,
\end{align*}\end{linenomath}
\rev{corresponding to the location of the tunnel in the depth 300 m below Earth's surface.
The initial pressure in the domain before excavation is given by the steady-state solution of the flow problem, satisfying the same boundary conditions as above and $p_D=300$ on $\Gamma^t$.}
The principal values of the initial stress $\vc\sigma_0$ are 60 MPa in the direction of excavation, 45 MPa in the orthogonal horizontal direction, and 11 MPa in the vertical direction.
For the values of other physical parameters of the model, we refer to \cref{tab:tunnel-params}.
\rev{The parameter values were inspired by the Tunnel Sealing Experiment (TSX) \cite{rutqvist_modeling_2009}.}
\begin{table}[h]
    \centering
    \caption{Physical parameters of the \rev{tunnel} excavation test problem.}
    \label{tab:tunnel-params}
    \begin{tabular}{l|c|r|c}
    \multicolumn{2}{c}{Quantity} & \multicolumn{1}{c}{Value} & Unit\\
    \hline
    hydraulic conductivity of intact rock & $k_0$ & $3\cdot10^{-13}$ & m/s\\
    hydraulic conductivity near excavation & $k_1$ & $10^{-8}$ & m/s\\
    storativity & $S_m,S_f$ & $7.29\cdot10^{-8}$ & 1/\rev{m}\\
    rock Young modulus & $E_m$ & $60\cdot10^9$ & Pa\\
    fracture Young modulus & $E_f$ & $60\cdot10^6$ & Pa\\
    Poisson's ratio & $\nu_m,\nu_f$ & 0.2 & $-$\\
    roughness coefficient & $\eta$ & 0.01 & $-$\\
    minimal \rev{fracture width} & $\delta_{min}$ & $10^{-6}$ & m\\
    Biot-Willis coefficient & $\alpha_m,\alpha_f$ & $0.2$ & $-$\\
    fluid density & $\varrho$ & $10^3$ & kg/m${}^3$\\
    fluid viscosity & $\mu$ & $10^{-3}$ & Pa$\cdot$s\\
    \end{tabular}
\end{table}

The simulation involves two stages: 40 days of excavation with the speed $v^e=1$ m/day and 320 days of relaxation.
\rev{For the whole simulation, 71 time steps were used, ranging from $\Delta t=1$ d in the first 40 days to $\Delta t=15$ d from day 120 till the end.}
The excavation process is simulated by setting time-dependent boundary conditions on the \rev{tunnel} surface:
\begin{linenomath}\begin{align*}
p_D(t,\vc x) &=300, & \vc t_N(t,\vc x) &= \vc\sigma_0\nn && \mbox{if }\vc x\in\Gamma^{\rev{t}},~x_1\ge v^e t,\\
p_D(t,\vc x) &=0, & \vc t_N(t,\vc x) &= \vc0 && \mbox{if }\vc x\in\Gamma^{\rev{t}},~x_1<v^e t.
\end{align*}\end{linenomath}
In cases where a fracture crosses the boundary of the 3D domain, we prescribe the same boundary condition on the fracture boundary as on the adjacent part of $\partial\Omega_m$.
On interior fracture boundaries, we set $\vc t_N=\vc0$.

The hydraulic conductivity $\tn K_m$ is considered isotropic in the whole domain.
For the intact rock, we use the value $k_0=3\cdot 10^{-13}$ m/s. In the zone influenced by the excavation, i.e. for $\vc x\in\Omega_m$ s.t. $0<x_1<\min\{v^et,40\}$, the value is increased up to the near-surface value $k_1=10^{-8}$ m/s depending on the distance $r$ from the \rev{tunnel} axis:
\[ k^e(r) = \begin{cases} k_1 & \mbox{if }r\le R+L_{in},\\ c_1 e^{c_2 r} & \mbox{if }R+L_{in}<r<R+L_{out},\\k_0 & \mbox{if } r\ge R+L_{out}.\end{cases} \]
Here $R=1.1$ m is the \rev{tunnel} radius, $L_{in}=1$ m, and $L_{out}=4$ m determine the interval of transition from $k_1$ to $k_0$. The parameters $c_1,c_2$ are chosen so that $k^e$ is continuous.

\begin{figure}
    \centering
    \includegraphics[width=0.49\textwidth]{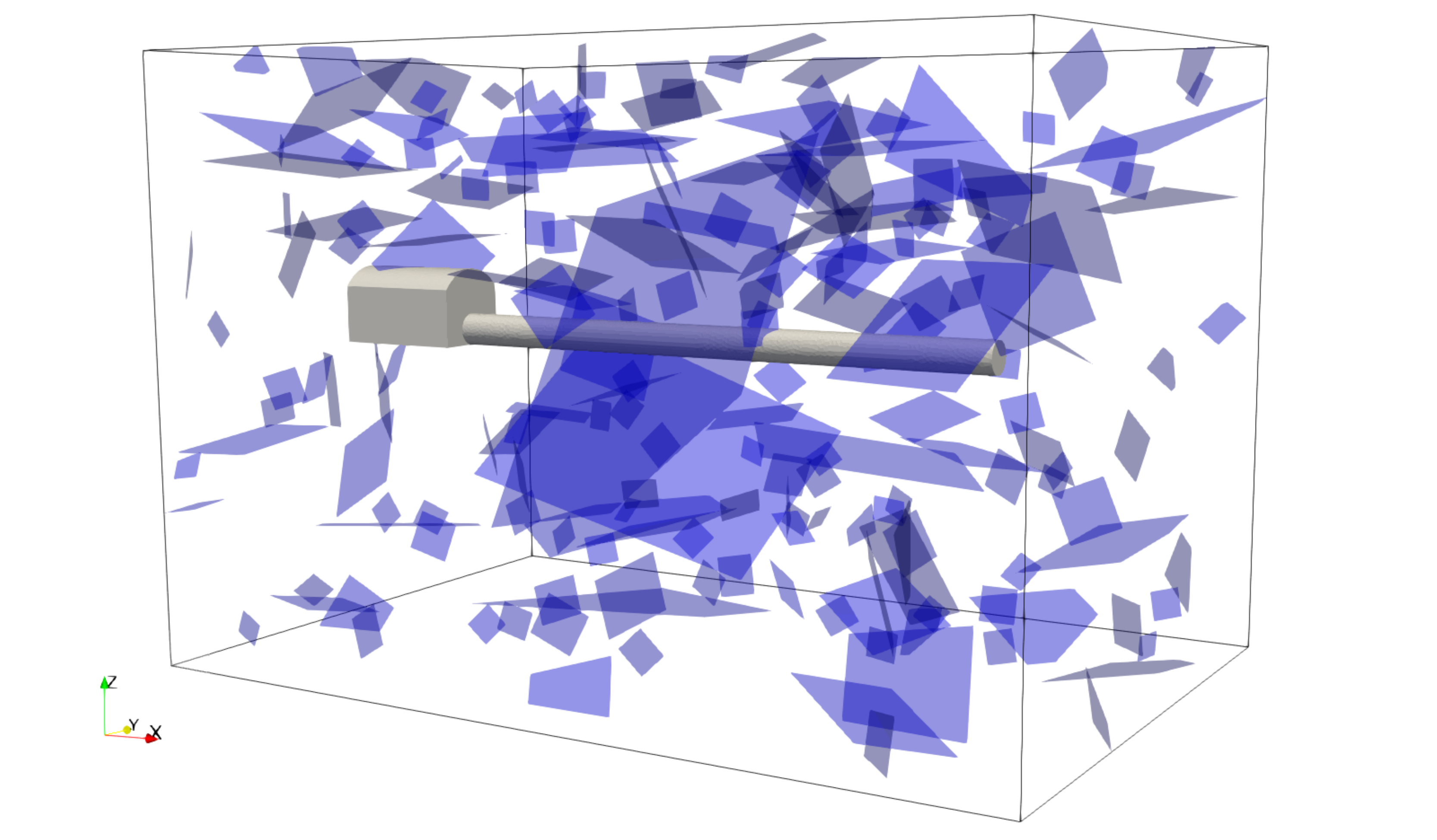}
    \includegraphics[width=0.49\textwidth]{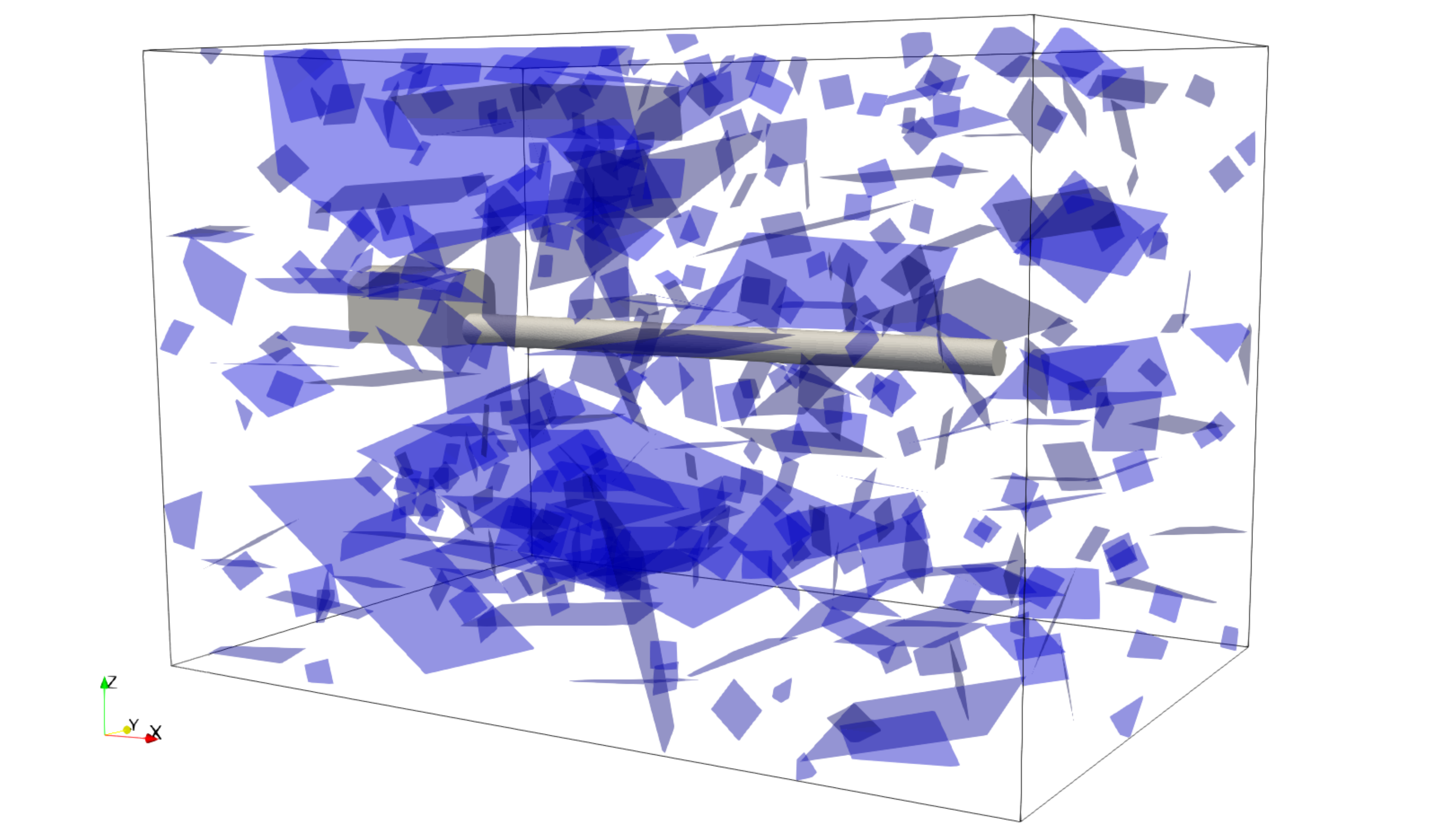}
    \caption{Two DFN configurations in the \rev{tunnel} excavation test. Left: 200 fractures, right: 400 fractures.}
    \label{fig:tunnel-geo-frac}
\end{figure}

\begin{figure}
    \centering
    \includegraphics[width=\textwidth]{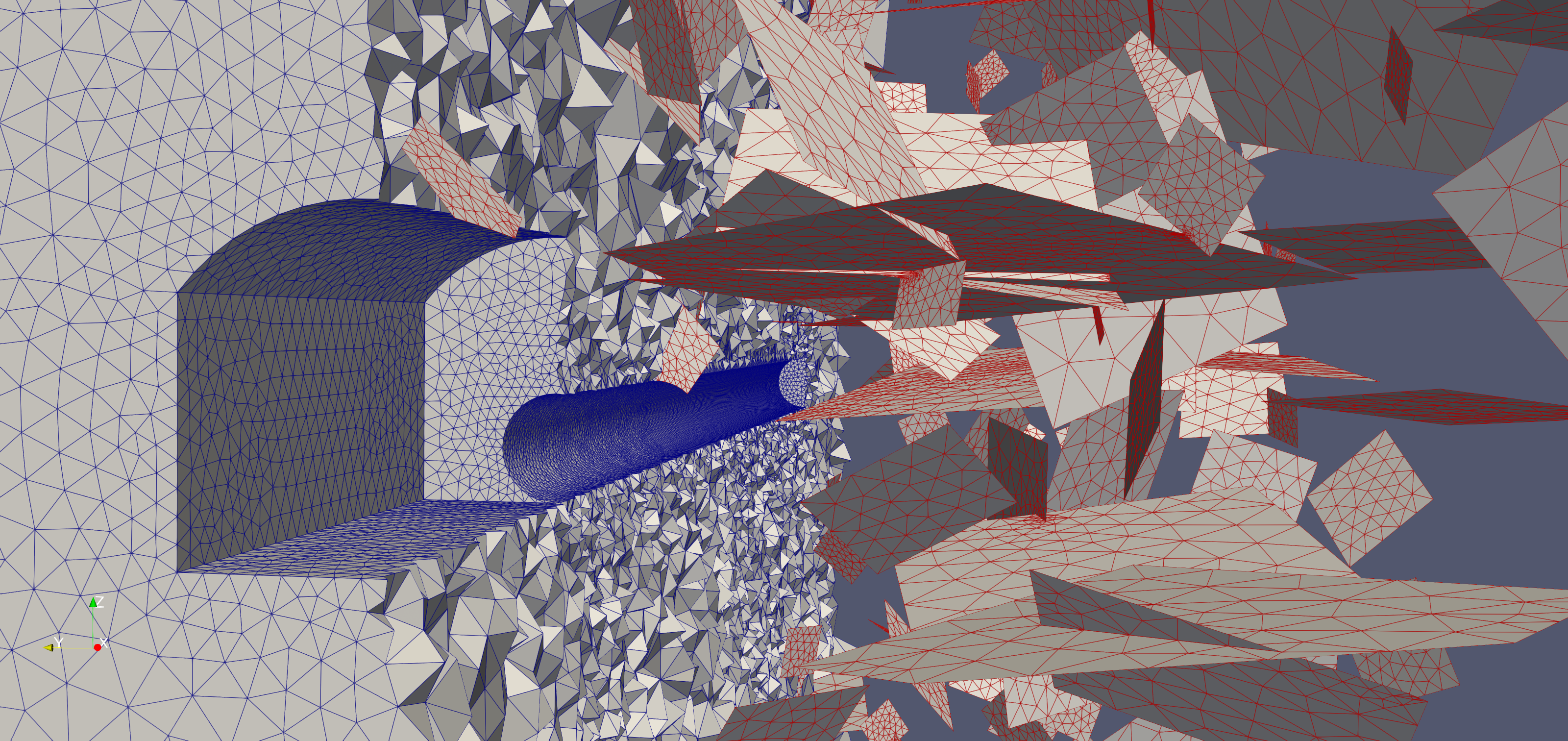}
    \caption{A detailed view of the computational mesh of the domain with 400 fractures. }
    \label{fig:mesh_400}
\end{figure}

The stochastic DFN model was taken from \cite{Ohman_Site_2010}, data for HDR -- repository domain, locality Forsmark, Sweden.
Two sets of square-shaped fractures were generated for the target set sizes 200 and 400, see \cref{fig:tunnel-geo-frac}.
The cross-sections $\delta$ vary between $0.1-0.6$ mm depending on the fracture size.
The initial hydraulic conductivity of the fractures $k_f = (\eta \varrho g/ 12\mu) \delta^2$ is given by the cross-section $\delta$ through the cubic law with the roughness parameter $\eta=0.01$. 
The roughness parameter was chosen to obtain realistic cross-sections for prescribed fracture transmissivities $T_f = k_f \delta$ fitted to the hydraulic tests in \cite{Ohman_Site_2010}.
The computational meshes for both DFN sets consist of approximately one million tetrahedral elements (see \cref{fig:mesh_400}).

\begin{figure}
    \centering
    \begin{subfigure}[b]{\textwidth}
    \includegraphics[width=0.49\textwidth]{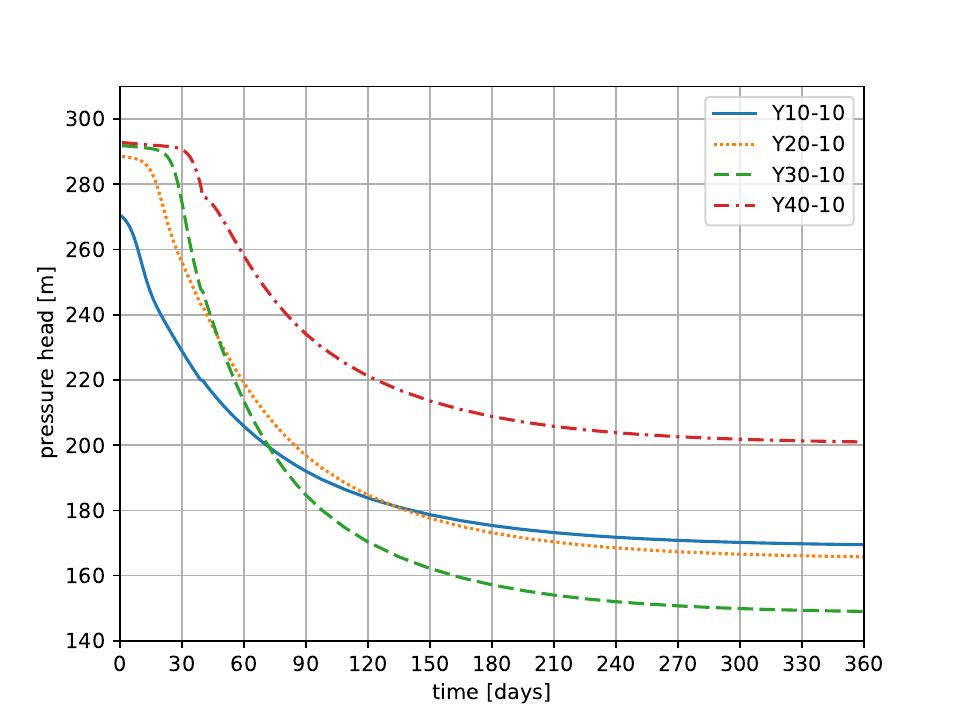}
    \includegraphics[width=0.49\textwidth]{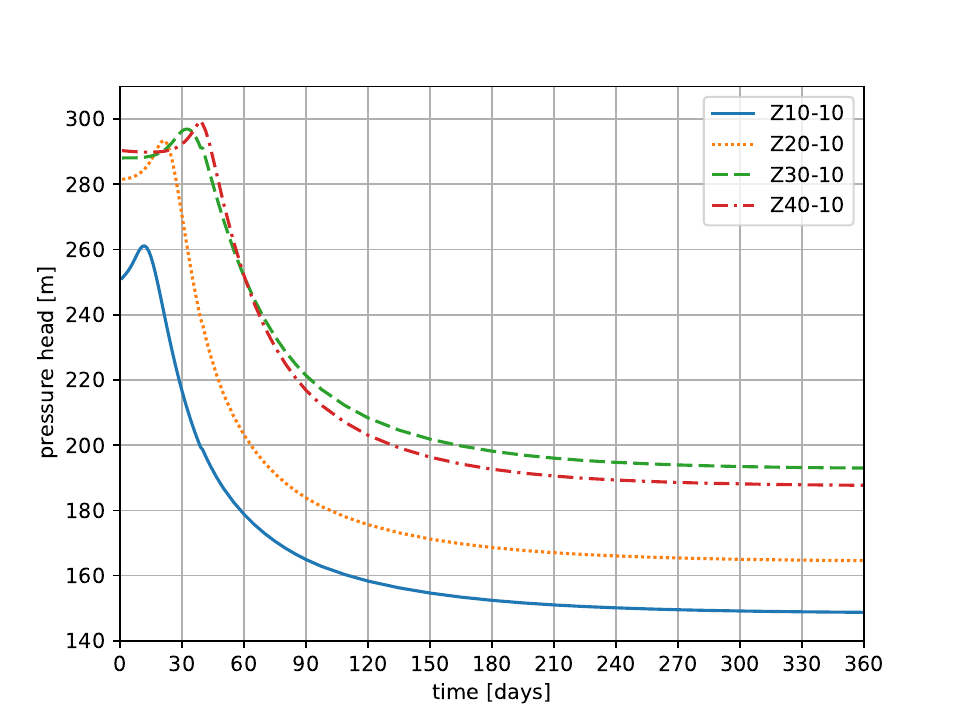}
    \caption{DFN with 200 fractures.}
    \end{subfigure}
    \begin{subfigure}[b]{\textwidth}
    \includegraphics[width=0.49\textwidth]{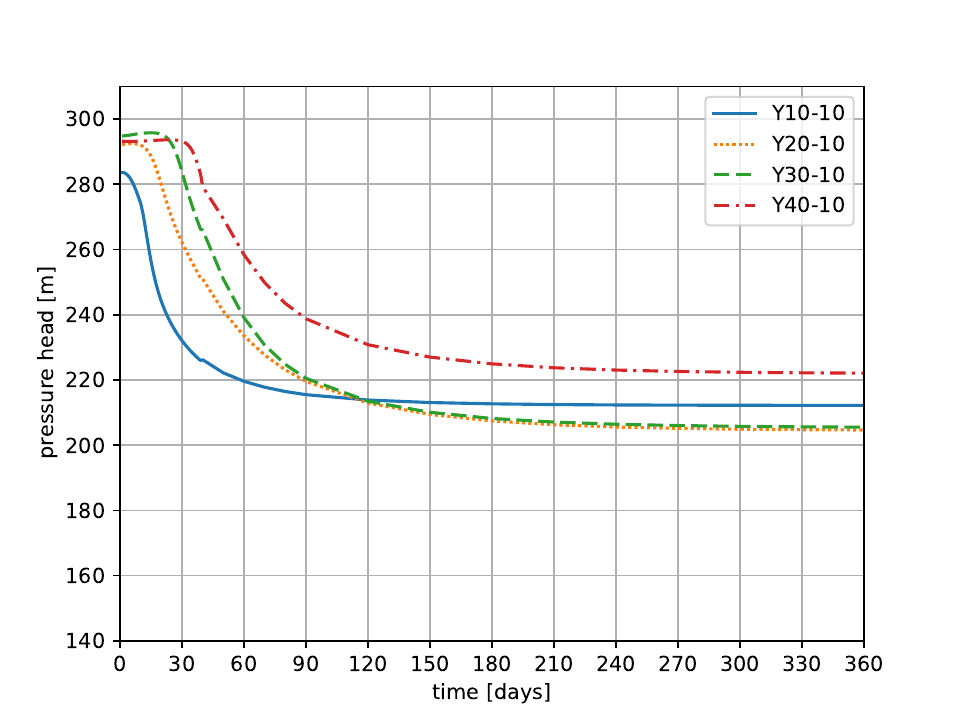}
    \includegraphics[width=0.49\textwidth]{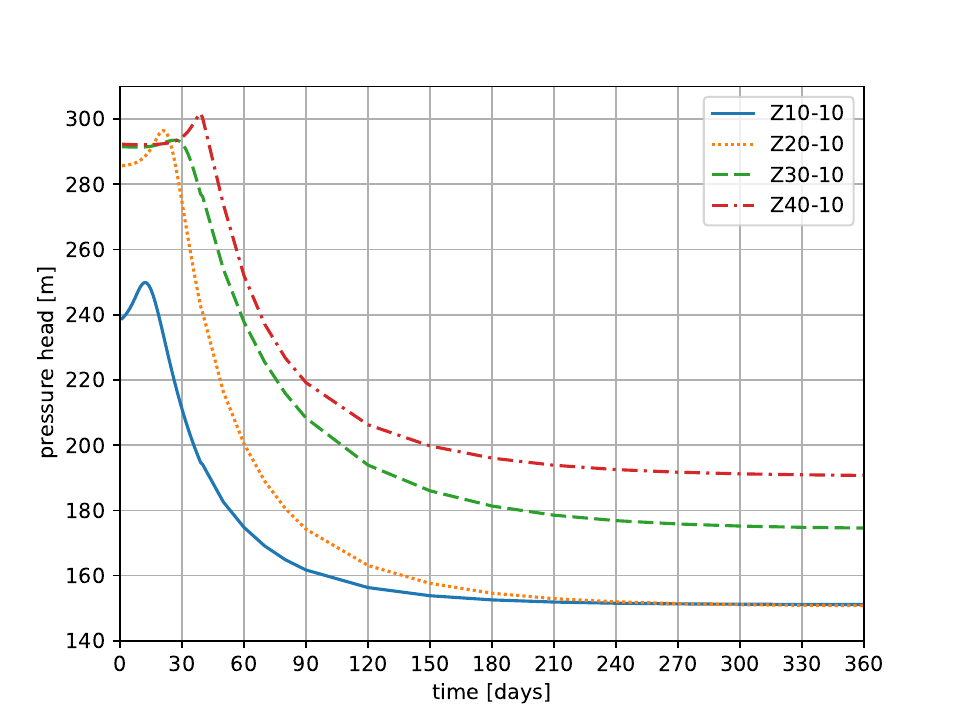}
    \caption{DFN with 400 fractures.}
    \end{subfigure}
    \begin{subfigure}[b]{\textwidth}
    \includegraphics[width=0.4\textwidth]{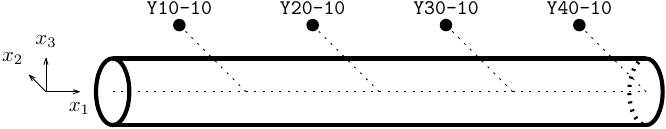}
    \hspace{0.1\textwidth}
    \includegraphics[width=0.4\textwidth]{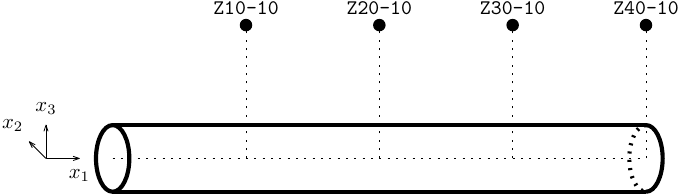}
    \caption{Location of observation points.}
    \end{subfigure}
    \caption{Evolution of pressure head in observation points located 10 m from the \rev{tunnel}. Left: observation points on the side of the \rev{tunnel}; Right: observation points above the \rev{tunnel}.}
    \label{fig:tunnel_pressure_obs}
\end{figure}

\begin{figure}
    \centering
    \begin{subfigure}[b]{0.49\textwidth}
    \includegraphics[width=\textwidth]{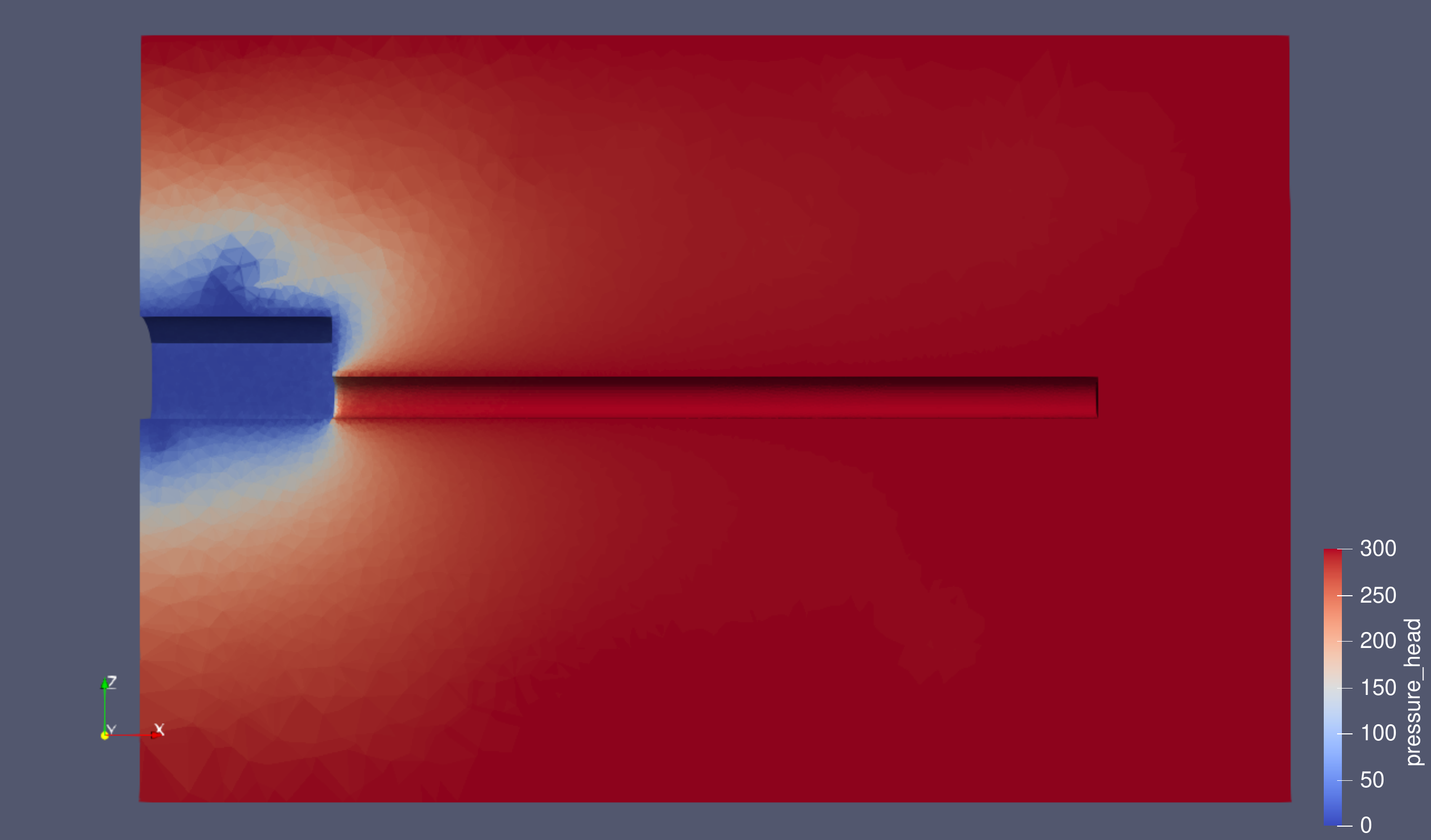}
    \includegraphics[width=\textwidth]{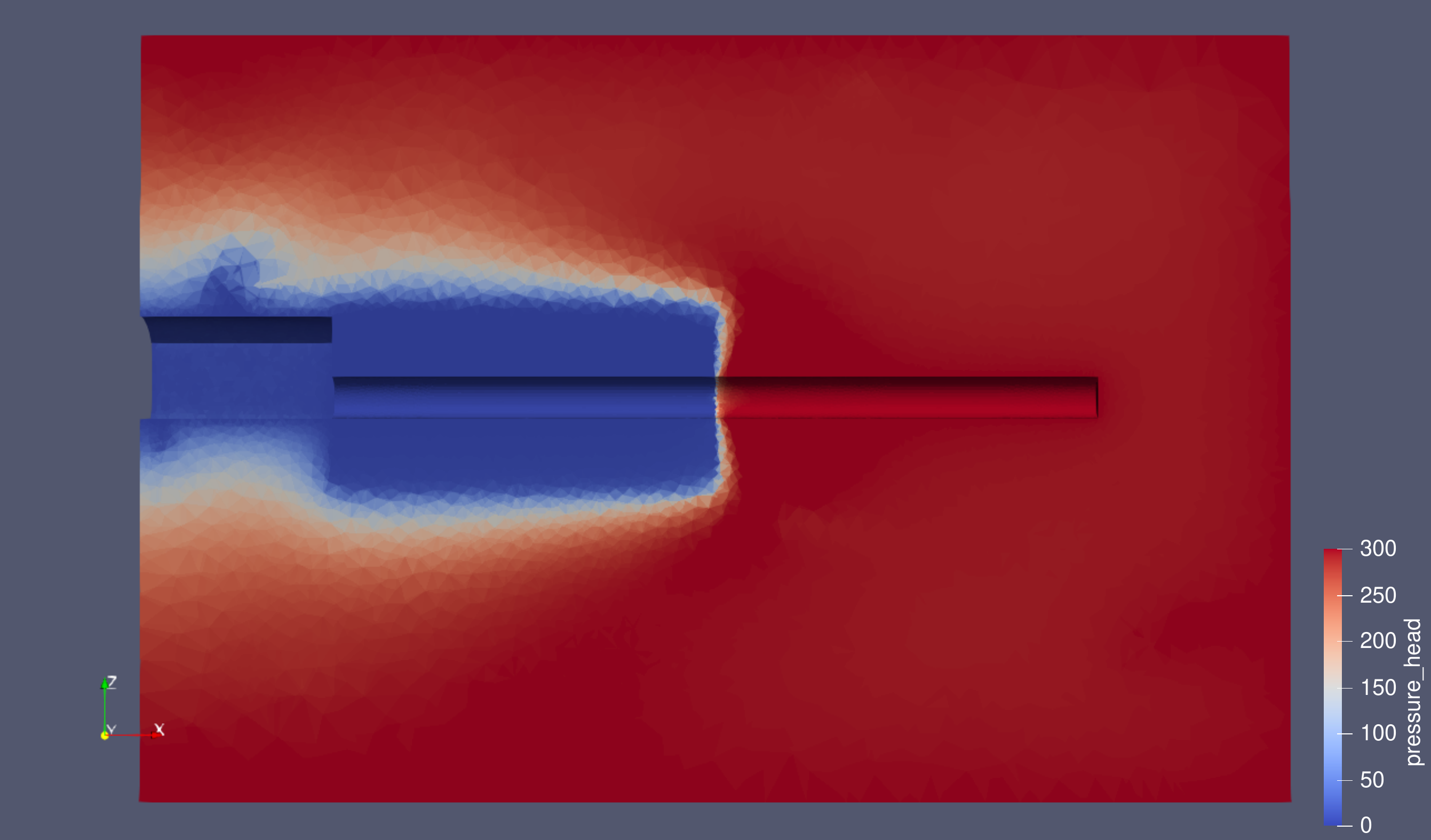}
    \includegraphics[width=\textwidth]{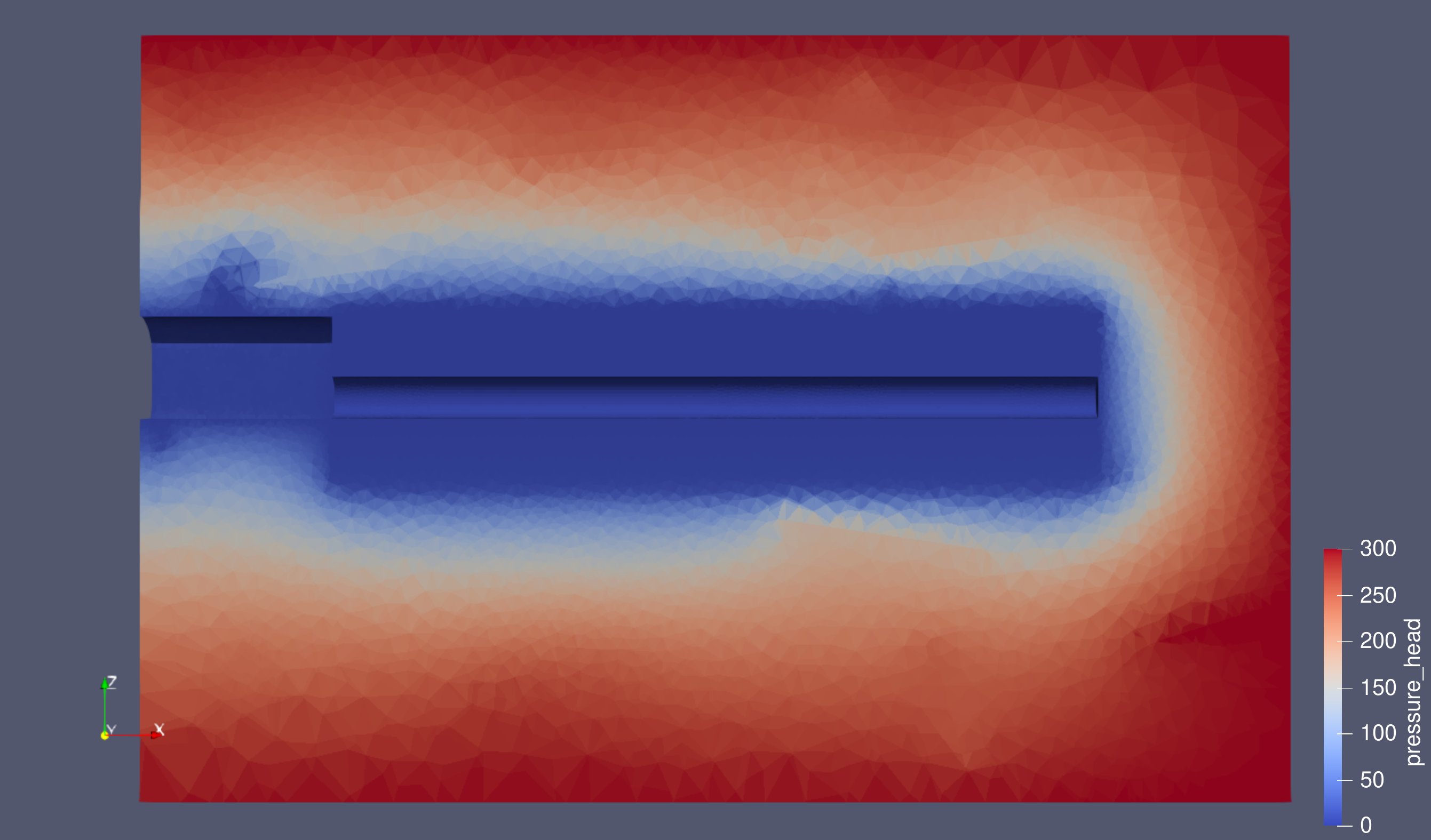}
    \caption{DFN with 200 fractures.}
    \end{subfigure}
    \begin{subfigure}[b]{0.49\textwidth}
    \includegraphics[width=\textwidth]{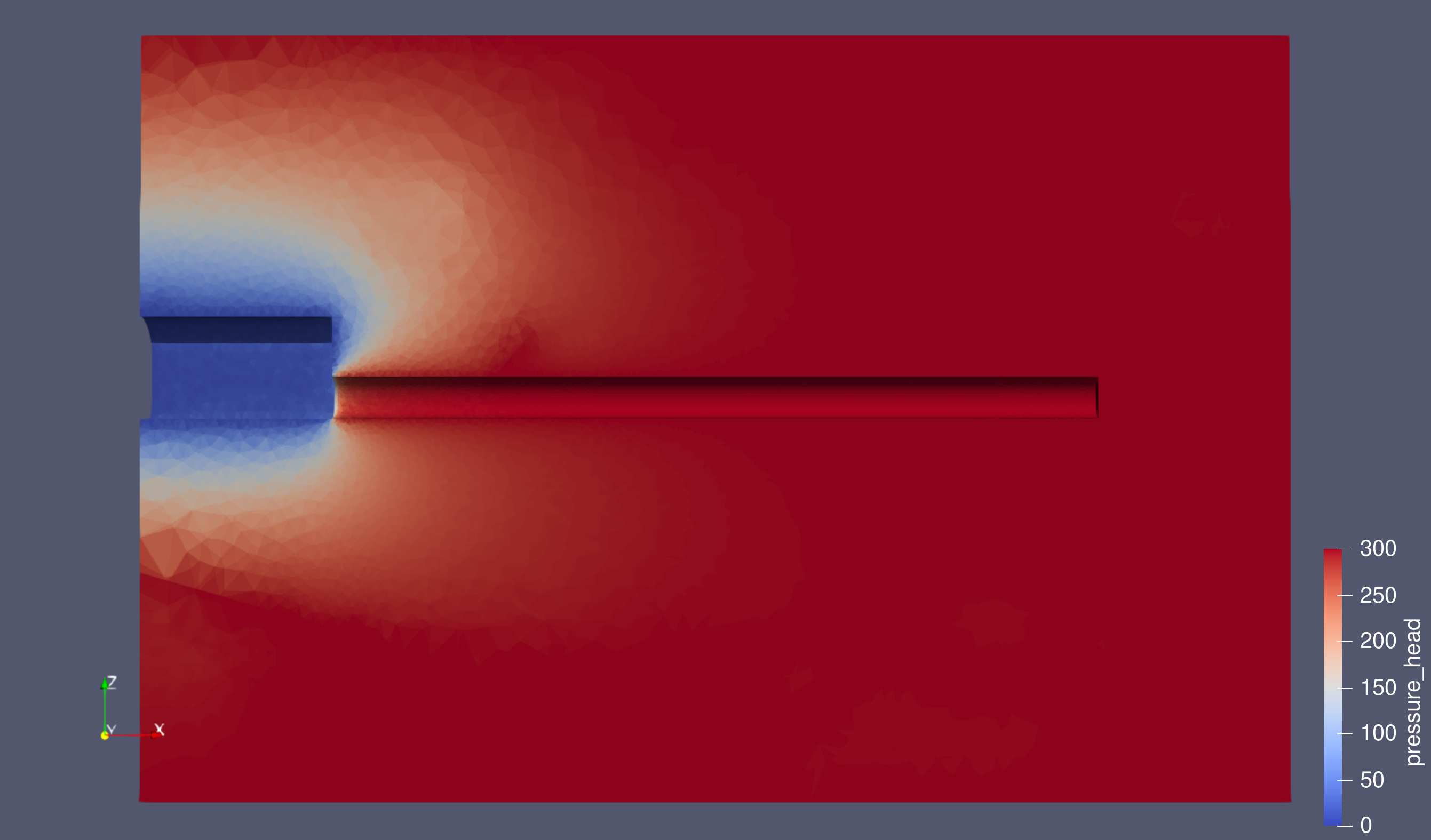}
    \includegraphics[width=\textwidth]{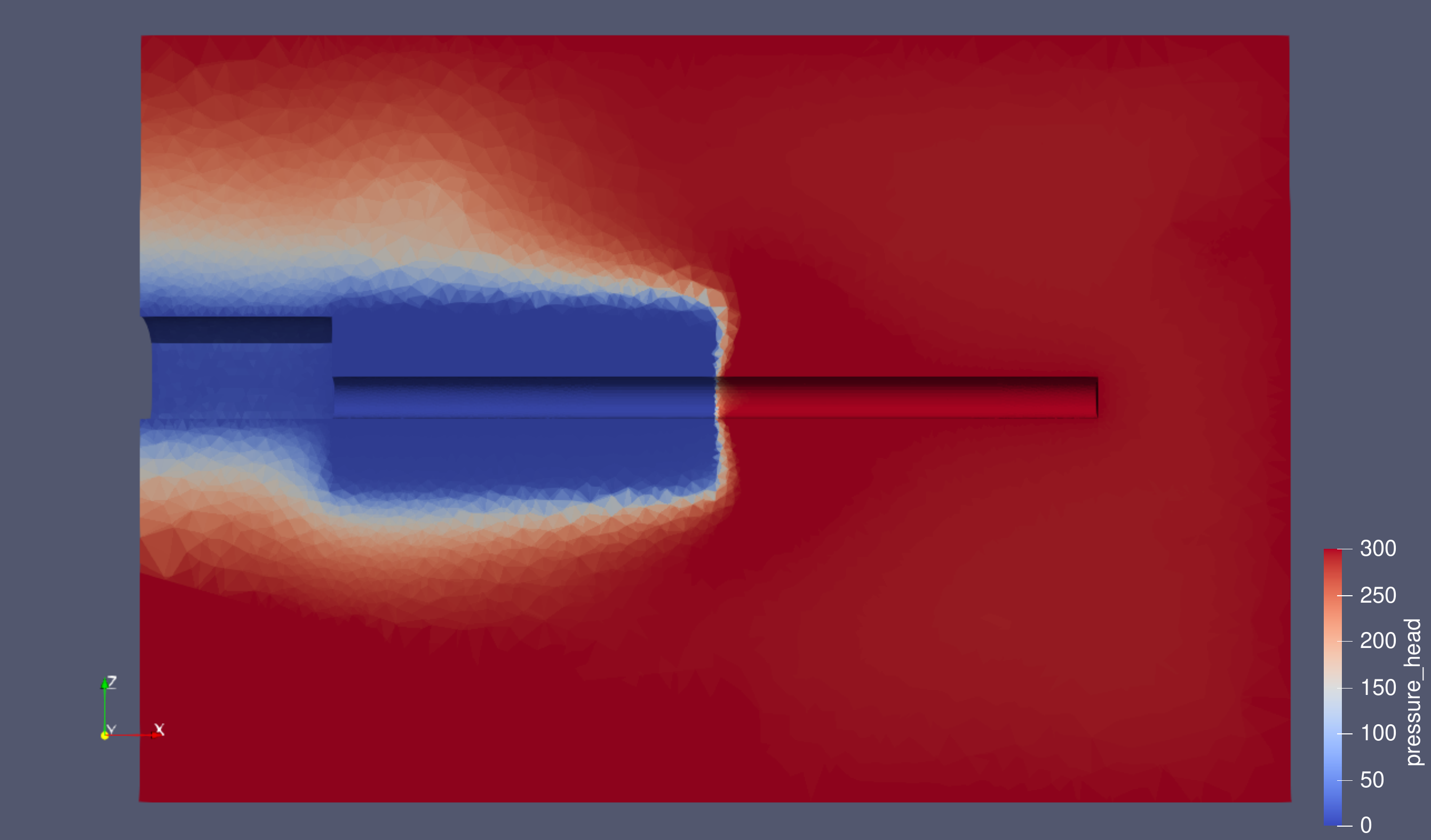}
    \includegraphics[width=\textwidth]{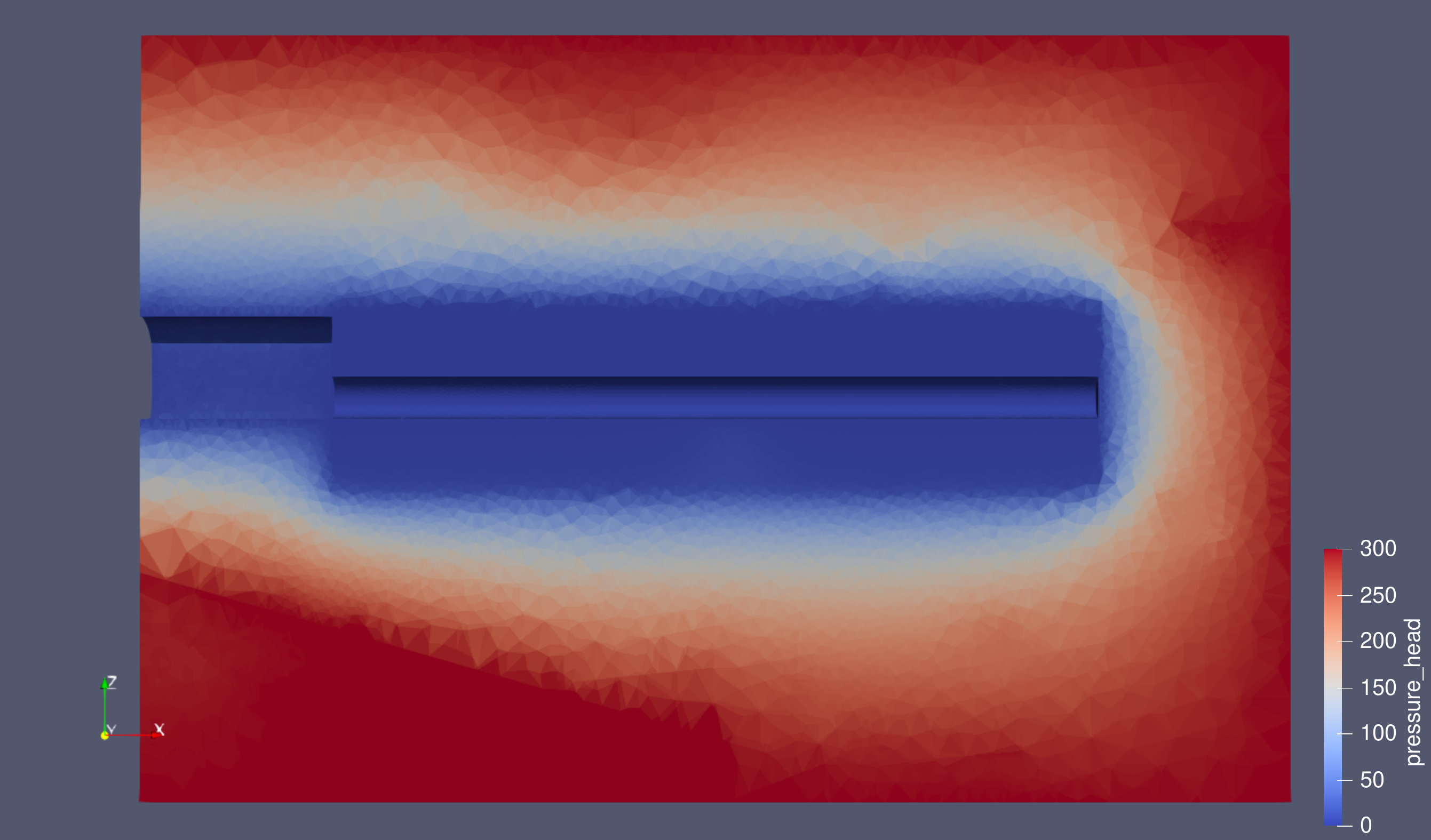}
    \caption{DFN with 400 fractures.}
    \end{subfigure}
    \caption{Comparison of pressure head distribution. From top to bottom: initial time, after 20 days, and after 360 days.}
    \label{fig:tunnel_pressure_plot}
\end{figure}

\begin{figure}
    \centering
    \begin{subfigure}[b]{\textwidth}
    \includegraphics[width=\textwidth]{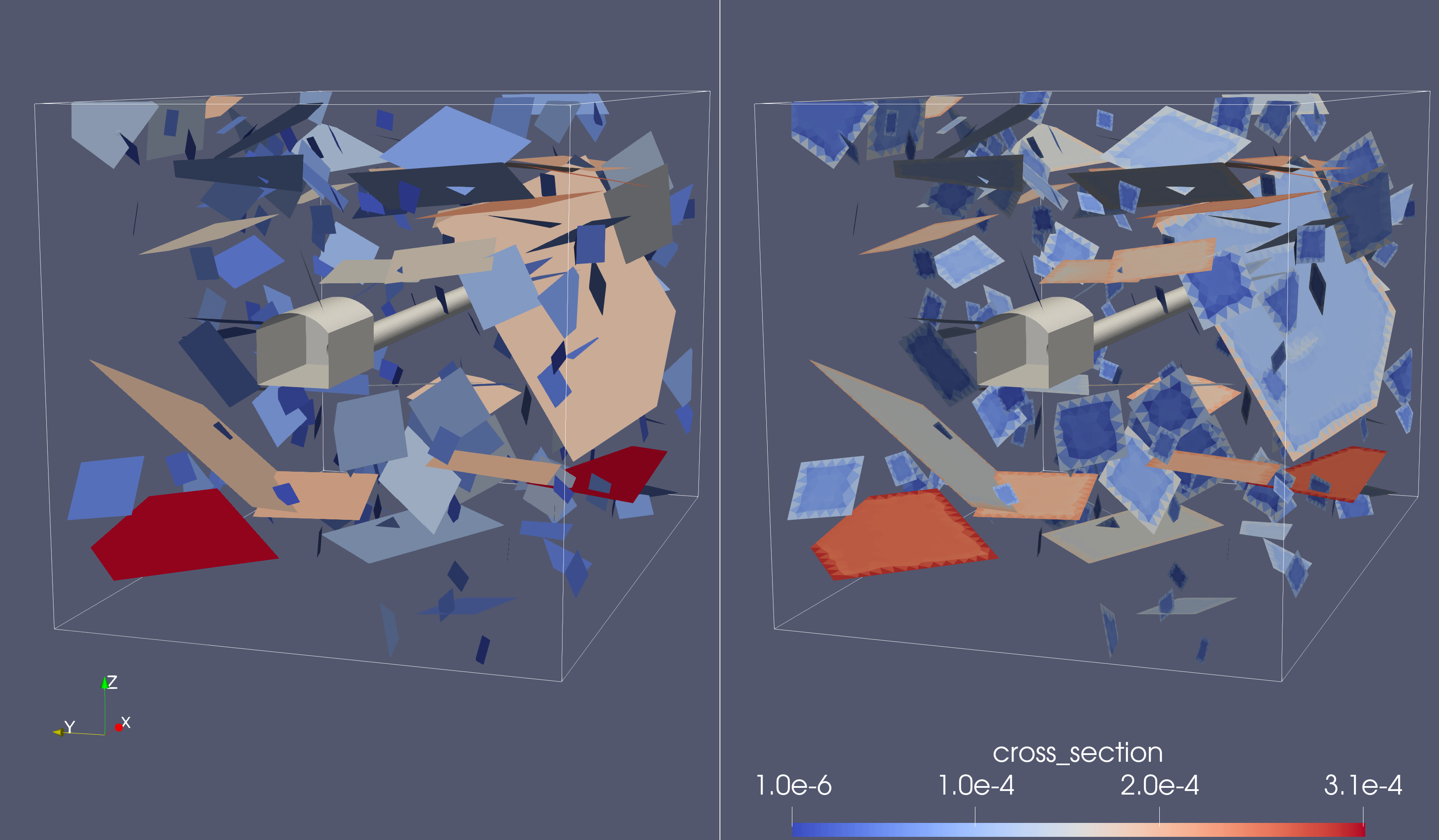}
    \caption{Initial cross-sections.\qquad\qquad\qquad\qquad (b) Cross-sections after 20 days.}
    \end{subfigure}
    \caption{DFN with 200 fractures: Comparison of initial fracture cross-sections (left) and cross-sections after 20 days of excavation (right).}
    \label{fig:tunnel-hcond}
\end{figure}

In \cref{fig:tunnel_pressure_obs}, the time evolution of the pressure head in several observation points is depicted.
The pressure above the \rev{tunnel} first increases until the excavation reaches the observation point and then gradually decreases. On the side of the \rev{tunnel}, the initial increase of the pressure is insignificant. This is in agreement with experimental observations e.g. from TSX \cite{rutqvist_modeling_2009}. Comparing the test case with 200 and 400 fractures, one can clearly see the difference in the final pressure head in the observation points.
The pressure field in the axial cut of the excavation is depicted in \cref{fig:tunnel_pressure_plot}.
While the overall behaviour of the pressure head is similar for both cases \rev{(200 and 400 fractures)}, one can observe differences in several parts of the domain, which are due to the high permeability of adjacent fractures.
\rev{The instant pressure drop is observed in the highly permeable zone around the excavation. At the final time, a smooth pressure distribution is visible as a result of the long-time relaxation. At this stage, the differences between the two cases become more evident.}
The comparison of initial fracture cross-sections and cross-sections during excavation can be seen in \cref{fig:tunnel-hcond}.
The initial stress induces the closing of most of the fractures to the minimal value $\delta_{min}=1~\mu$m and, consequently, the decrease of their hydraulic conductivity.
In summary, the effect of mechanics as well as of fracture network cannot be neglected, and it is important to consider the coupled model.

\begin{figure}
    \centering
    \includegraphics[width=0.49\textwidth]{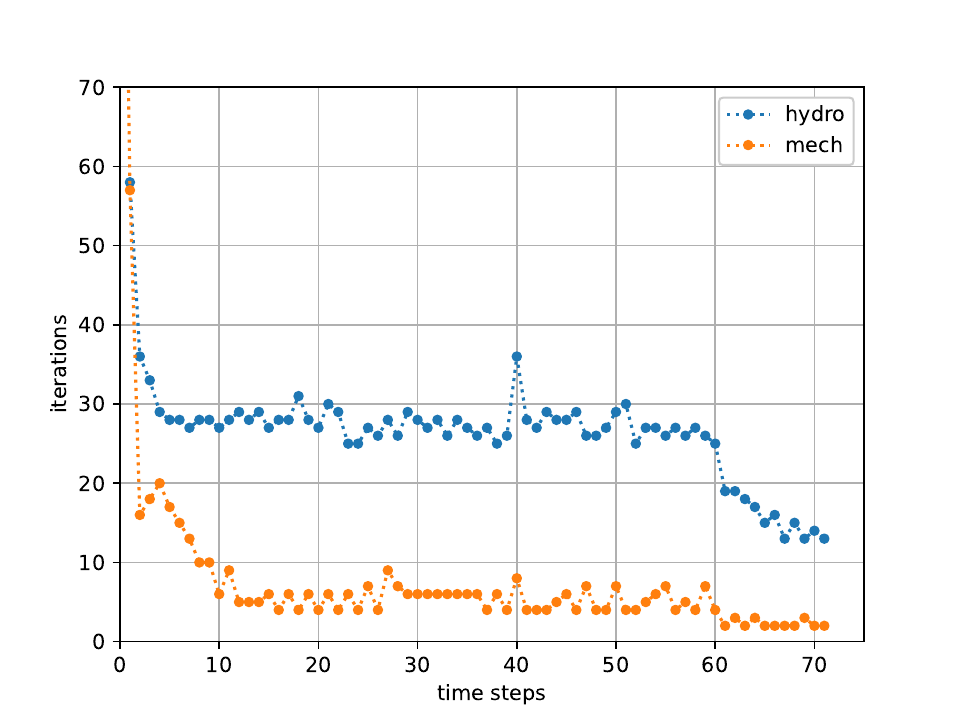}
    \includegraphics[width=0.49\textwidth]{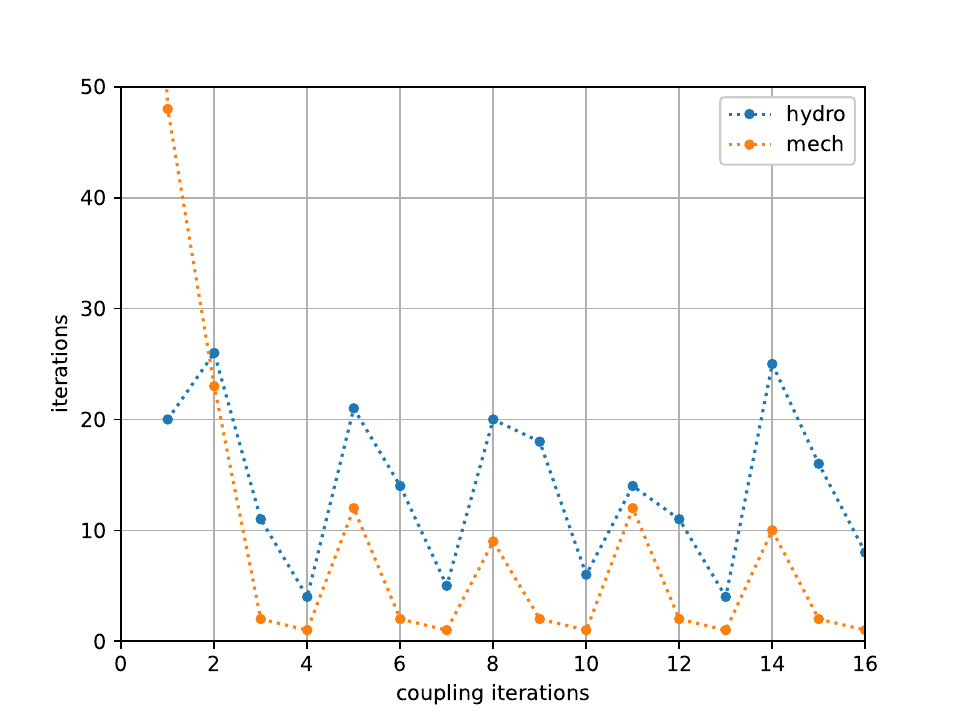}
    \caption{Cumulative number of solvers iterations in each time step (left) and the number of solvers iterations for each coupling iteration in the first 5 time steps (right) \rev{for the problem with one million elements and 200 fractures}. Time step/coupling iteration zero represents the steady state initial solution, which took 148 iterations.}
    \label{fig:iters}
\end{figure}

We also investigated the numerical efficiency of our solvers on the meshes described above as well as varying the mesh size for the model with 200 fractures.
The simulations were run on 1 node of the \rev{LUMI} \cite{lumi} supercomputer using \rev{64} cores. \rev{The relative tolerances were set to  $10^{-9}$ and $10^{-6}$ for the hydrological and the mechanical subproblem solvers, respectively, and $10^{-4}$ for the iterative coupling.} The first time step needed 4 iterations of the hydro-mechanical coupling, while the rest of the simulation needed 3 iterations except for the last few time steps going down to 2. Overall, the computations required 201-203 iterations of the hydro-mechanical coupling for the 71 time steps. The cumulative number of iterations needed by the hydrological subproblem solver, the number of Hessian multiplications done by the mechanical subproblem solver, and the approximate simulation time are reported in \cref{tab:num_scal}. We report the number of the Hessian multiplications for MPGP instead of the number of iterations due to it being the most expensive part of the algorithm that is carried out either once (by CG or proportioning step) or twice (by expansion step) per iteration.
Both solvers employ a warm starting strategy that uses the previous solution  (dual solution in the case of the mechanical subproblem) as the initial guess. \cref{fig:iters} shows the cumulative number of solvers iterations in each time step (left) and illustrates the effect of warm starting on the first 5 time steps (right) for the problem with about one million elements and 200 fractures. \rev{As expected, the number of iterations for the mechanical subproblem grows modestly with the number of elements in the mesh. However, the algorithm seems to be numerically scalable with respect to the number of fractures.}

\begin{table}
    \centering
    \caption{Numerical scalability for each test case. The problem name consists of the \rev{number of elements in thousands and the number of fractures in the brackets}. We report the primal and the dual dimensions of the mechanical subproblem, the number of hydraulic subproblem iterations, the number of Hessian multiplications needed by the mechanical subproblem solver, and the approximate run time in minutes. \rev{In all cases, 71 time steps were performed.}}
    \label{tab:num_scal}
    \begin{tabular}{l|r|r|r|r|r}
    Problem & Primal & Dual & Hydr. iter. & Mech. iter. & Time [min]\\
    \hline
    \rev{306k (200)} & 213,945 & 18,632 & \rev{2,108} & \rev{365} & \rev{6.4}\\
    \rev{502k (200)} & 322,953 & 19,591 & \rev{2,587} & \rev{544} & \rev{10.3}\\
    \rev{1052k (200)} & 629,352 & 28,128 & \rev{1,869} & \rev{625} & \rev{20.2}\\
    \rev{989k (400)} & 651,861 & 42,759 & \rev{1,894} & \rev{570} & \rev{19.7}
    \end{tabular}
\end{table}

\section{Conclusions}\label{sec:conclusions}
We presented a numerical method for the solution of flow and mechanics in fractured rock, which combines the continual description of rock mass with the discrete fracture network.
\rev{
The method was implemented using a suitable combination of state-of-the-art open-source libraries and in-house codes with the aim of being capable of solving problems with a high number of fractures.
From the performed computational tests we make the following conclusions and remarks.
\begin{itemize}
\item \rev{The proposed nonlinear hydro-mechanical model of fractured rock has been validated against a 2D single-fracture nonlinear problem that admits a semi-analytical solution.}
\item The finite element discretization yields satisfactory results for a 3D benchmark with realistic physical parameters. However, for problems in a nearly incompressible regime, further improvement may be necessary to avoid pressure oscillations, such as the algebraic stabilization proposed in \cite{frigo2021efficient}.
\item \rev{The robust implementation of contact mechanics with convergence efficiency comparable to the solution of the hydraulic subproblem marks a significant advancement towards a more realistic model by capturing the non-penetration condition. A logical next step is the integration of friction laws, which can be included by an appropriate modification of \cref{eq:contact_friction}. However, the iterative scheme for the coupled problem must be revisited, and the related quadratic programming problems may require more advanced modifications of the MPGP algorithm.}
\item \rev{
The convergence of the hydro-mechanical coupling proved to be robust with respect to mesh size. However, using a deformation-dependent fracture aperture adversely affected the convergence of the iterative splitting and the conditioning of the algebraic problems. This issue became apparent in the 3D benchmark, leading to the adoption of a fixed fracture aperture in this case. Addressing these challenges will necessitate further investigation and detailed analysis.}
\item The proposed computational approach can solve problems with \rev{millions of elements and} hundreds of fractures. To \rev{substatially increase the number of elements}, further improvement of the parallel scalability of the method is necessary. For this reason, a suitable domain decomposition technique, such as FETI \cite{Farhat1991,DosHorKuc-CNME-06}, can be considered.
\end{itemize}
}


\section*{CRediT authorship contribution statement}

{\bf Jan Stebel:} Methodology, Software, Investigation, Visualization, Writing - original draft,
{\bf Jakub Kružík:} Software, Investigation, Visualization, Writing - original draft,
{\bf David Horák:} Software, Writing - review  editing,
{\bf Jan Březina:} Software, Writing - review \& editing,
{\bf Michal Béreš:} Validation, Writing - review \& editing.

\section*{Data availability}
Data will be made available on request.

\section*{Declaration of Competing Interest}
The authors declare that they have no known competing financial interests or personal relationships that could have appeared to influence the work reported in this paper.

\section*{Acknowledgement}
The research was supported by European Union’s Horizon 2020 research and innovation programme under grant agreement number 847593 and by The Czech
Radioactive Waste Repository Authority (SÚRAO) under grant agreement
number SO2020-017. 

The work of J. Kru\v{z}\'ik, D. Hor\'ak and M. B\'{e}re\v{s} is partially supported by the financial support of the European Union under the REFRESH - Research Excellence For REgion Sustainability and High-tech Industries project number \url{CZ.10.03.01/00/22_003/0000048} via the Operational Programme Just Transition.

The work was partially supported by Grant of SGS No. SP2024/067, VŠB - Technical
University of Ostrava, Czech Republic.

Computational resources were provided by the e-INFRA CZ projects ID:90140 and ID:90254, supported by the Ministry of Education, Youth and Sports of the Czech Republic. 

We acknowledge IT4Innovations, Czech Republic, for awarding this pro\-ject access to the LUMI supercomputer, owned by the EuroHPC Joint Undertaking, hosted by CSC (Finland) and the LUMI consortium.

\appendix

\section{Weak formulation of DFM hydro-mechanical problem}
\label{sec:appendix_weak}
In this section, we present and derive the weak formulation of the Problem (DFM-HM).

\paragraph{Weak formulation}
A triplet $(\uu,p,\vv)\in L^2(I;\mathcal K\times\mathcal Q\times\mathcal W(v_N))$ is said to be the weak solution of Problem (DFM-HM) if $p(0,\cdot)=p_0$ in $\Omega_m\times\Omega_f$ and, for any $(\zz,q,\ww)\in\mathcal K\times\mathcal Q\times\mathcal W(0)$ and a.a. in $I$, the following is satisfied:
\begin{linenomath}\begin{align}
\label{eq:weak_var_ineq} a(\uu,\zz-\uu) - b(\zz-\uu,p) &\ge f_1(\zz-\uu),\\
\label{eq:weak_flow} \dt(c(p,q) + b(\uu,q)) + d(\vv,q) &= f_2(q),\\
\label{eq:weak_darcy} e(\uu;\vv,\ww) - d(\ww,p) &= f_3(\ww).
\end{align}\end{linenomath}
Here, the sets of functions are defined as follows:
\begin{linenomath}\begin{align*}
\mathcal K&:=\{(\zz_m,\zz_f)\in H^1(\Omega_m;\Real^d)\times H^1(\Omega_f;\Real^d);\\
&\qquad\qquad\zz_*=\uu_{D*}~\mbox{on}~\Gamma_{D*}^{mech},~*\in\{m,f\},~a_f(\zz)\ge\delta_{min}\},\\
\mathcal Q&:=L^2(\Omega_m)\times L^2(\Omega_f),\\
\mathcal W(\varphi)&:=\{(\ww_m,\ww_f)\in\Hdiv(\Omega_m)\times\Hdiv(\Omega_f);\\
&\qquad\qquad\ww_*\cdot\nn=\varphi_* \mbox{ on }\Gamma_{N*}^{flow},~*\in\{m,f\}\},
\end{align*}\end{linenomath}
where $\varphi=(\varphi_m,\varphi_f)$ is an arbitrary pair of functions defined in $\Gamma_{Nm}^{flow}\times\Gamma_{Nf}^{flow}$.
The forms used in the weak formulation are:
\begingroup
\allowdisplaybreaks
\begin{linenomath}\begin{subequations}\label{eq:weak_forms}
\begin{align}
a(\uu,\zz) =& \int_{\Omega_m} \vc C_m\eps(\uu_m):\eps(\zz_m) + \int_{\Omega_f}\delta\avg{\vc C_f\feps(\uu):\feps(\zz)}, \\
b(\uu,p) =& \int_{\Omega_m} \alpha_m \varrho g p_m\div\uu_m + \int_{\Omega_f} \delta \alpha_f \varrho g p_f\fdiv\uu, \\
c(p,q) =& \int_{\Omega_m} S_m p_m q_m + \int_{\Omega_f} \delta S_f p_f q_f, \\
d(\vv,q) =& \int_{\Omega_m} q_m\div\vv_m + \int_{\Omega_f} q_f\fdiv(\delta\vv_m,\vv_f), \\
e(\uu;\vv,\ww) =& \int_{\Omega_m} \tn K_m^{-1}\vv_m\cdot\ww_m + \int_{\Omega_f} \frac1{\delta k_f(\uu)}\vv_f\cdot\ww_f\notag\\
&+ \int_{\Omega_f}\frac\delta{k_f}\avg{(\vv_m\cdot\nnu)(\ww_m\cdot\nnu)}, \\
f_1(\zz) =& \int_{\Omega_m} \ff_m\cdot\zz_m + \int_{\Gamma_{Nm}^{mech}}\vc t_{Nm}\cdot\zz_m + \int_{\Omega_f} \delta\ff_f\cdot\zz_f + \int_{\Gamma_{Nf}^{mech}} \delta\vc t_{Nf}\cdot\zz_f\notag\\
       & - \int_{\Omega_m}\vc\sigma_{0m}:\eps(\zz_m) - \int_{\Omega_f}\delta\vc\sigma_{0f}:\feps(\zz),\\
f_2(q) =& \int_{\Omega_m} s_m q_m + \int_{\Omega_f} \delta s_f q_f,\\
f_3(\ww) =& \int_{\Gamma_{Dm}^{flow}}p_{Dm}\ww_m\cdot\nn + \int_{\Gamma_{Df}^{flow}}p_{Df}\ww_f\cdot\nn
 - \int_{\Omega_m}\vc g\cdot\ww_m\notag\\
 &- \int_{\Omega_f}\vc g_\tau\cdot\ww_f - \int_{\Omega_f}\delta\avg{(\vc g\cdot\nnu)(\ww_m\cdot\nnu)}.
\end{align}
\end{subequations}\end{linenomath}
\endgroup
In what follows, we derive the variational inequality \cref{eq:weak_var_ineq} and the identities \cref{eq:weak_flow,eq:weak_darcy}.
We start with the mechanical part.

\paragraph{Weak form of balance of momentum \cref{eq:weak_var_ineq}}
Let us multiply \cref{eq:el} by $\zz_m-\uu_m$, where $\zz=(\zz_m,\zz_f)\in\mathcal K$, and integrate over $\Omega_m$.
Using Green's theorem and the boundary conditions on $\Gamma_{Dm}^{mech}$ and $\Gamma_{Nm}^{mech}$, we obtain the identity:
\begin{linenomath}\ml{\label{eq:weak_el_m}
\int_{\Omega_m}\sigmapor_m:\eps(\zz_m-\uu_m) + \int_{\Omega_f}2\avg{\sigmapor_m\nnu\cdot(\zz_m-\uu_m)}\\
= \int_{\Omega_m}\ff_m\cdot(\zz_m-\uu_m) + \int_{\Gamma_{Nm}^{mech}}\vc t_{Nm}\cdot(\zz_m-\uu_m). 
}\end{linenomath}
Similarly, we multiply \cref{eq:mixed_dim_problem_el_frac} by $\zz_f-\uu_f$ and integrate over $\Omega_f$:
\eq{\label{eq:weak_el_f0} -\int_{\Omega_f} \fdiv(\delta\sigmapor)\cdot(\zz_f-\uu_f) = \int_{\Omega_f}\delta\ff_f\cdot(\zz_f-\uu_f). }
In this equation, we express $\fdiv(\delta\sigmapor)$ as follows:
\begin{linenomath}\mls{ \fdiv(\delta\sigmapor)=\div_\tau(\delta\sigmapor_f)+\delta\avg{\sigmapor\nnu} = \div_\tau(\delta\sigmapor_f) + (\sigmapor_m^+-\sigmapor_f)\nnu^+ + (\sigmapor_m^--\sigmapor_f)\nnu^-\\
= \div_\tau(\delta\sigmapor_f)+2\avg{\sigmapor_m\nnu}, }\end{linenomath}
apply Green's theorem and employ the boundary conditions on $\Gamma_{Df}^{mech}$ and $\Gamma_{Nf}^{mech}$.
Then \cref{eq:weak_el_f0} can be rewritten in the following form:
\begin{linenomath}\ml{\label{eq:weak_el_f}
\int_{\Omega_f}\left(\delta\sigmapor_f:\eps_\tau(\zz_f-\uu_f) - 2\avg{\sigmapor_m\nnu}\cdot(\zz_f-\uu_f) \right)\\
= \int_{\Omega_f}\delta\ff_f\cdot(\zz_f-\uu_f) + \int_{\Gamma_{Nf}^{mech}}\delta\vc t_{Nf}\cdot(\zz_f-\uu_f). 
}\end{linenomath}
Adding \cref{eq:weak_el_m} and \cref{eq:weak_el_f}, we obtain:
\begin{linenomath}\ml{\label{eq:weak_el_sum}
\int_{\Omega_m}\sigmapor_m:\eps(\zz_m-\uu_m)
+\int_{\Omega_f}\delta\sigmapor_f:\eps_\tau(\zz_f-\uu_f) + \int_{\Omega_f}2\avg{\sigmapor_m\nnu\cdot(\zz_m-\uu_m-(\zz_f-\uu_f)) }\\
= \int_{\Omega_m}\ff_m\cdot(\zz_m-\uu_m) + \int_{\Gamma_{Nm}^{mech}}\vc t_{Nm}\cdot(\zz_m-\uu_m)
+\int_{\Omega_f}\delta\ff_f\cdot(\zz_f-\uu_f) + \int_{\Gamma_{Nf}^{mech}}\delta\vc t_{Nf}\cdot(\zz_f-\uu_f).
}\end{linenomath}
Now we will apply the contact conditions \crefrange{eq:contact_balance}{eq:contact_complementarity}.
Using the definition of $\vc\Lambda^\pm$, we can express the term under the third integral in \cref{eq:weak_el_sum} as follows:
\begin{linenomath}\ml{\label{eq:weak_el_term_sigma}
2\avg{\sigmapor_m\nnu\cdot(\zz_m-\uu_m-(\zz_f-\uu_f)) }\\
= 2\avg{\sigmapor_f\nnu\cdot(\zz_m-\uu_m-(\zz_f-\uu_f)) } - 2\avg{\vc\Lambda\cdot(\zz_m-\uu_m-(\zz_f-\uu_f)) }\\
= \delta\avg{\sigmapor_f:(\Delta(\zz-\uu)\otimes\nnu) } - 2\avg{\vc\Lambda\cdot(\zz_m-\uu_m)} + 2\underbrace{\avg{\vc\Lambda}}_{=0}\cdot(\zz_f-\uu_f),
}\end{linenomath}
where the last term vanishes due to \cref{eq:contact_balance}.
We further rewrite the last but one term in \cref{eq:weak_el_term_sigma}:
\begin{linenomath}\ml{\label{eq:weak_el_term_lambda}
2\avg{\vc\Lambda\cdot(\zz_m-\uu_m)}
\stackrel{\eqref{eq:contact_friction}}{=} 2\avg{(\vc\Lambda\cdot\nnu)((\zz_m-\uu_m)\cdot\nnu)}\\
\stackrel{\eqref{eq:contact_balance}}{=} 2\avg{\vc\Lambda\cdot\nnu}\avg{(\zz_m-\uu_m)\cdot\nnu}
= \avg{\vc\Lambda\cdot\nnu}(a_f(\zz)-a_f(\uu))\\
\stackrel{\eqref{eq:contact_complementarity}}{=} \avg{\vc\Lambda\cdot\nnu}(a_f(\zz)-\delta_{min})
\ge 0.
}\end{linenomath}
The last inequality follows from \cref{eq:contact_stress} and the fact that $\zz\in\mathcal K$.
Altogether, \crefrange{eq:weak_el_sum}{eq:weak_el_term_lambda} yield:
\begin{linenomath}\mls{
\int_{\Omega_m}\sigmapor_m:\eps(\zz_m-\uu_m) + \int_{\Omega_f}\delta\avg{\sigmapor_f:\feps(\zz-\uu)}
\ge \int_{\Omega_m}\ff_m\cdot(\zz_m-\uu_m)\\
+ \int_{\Gamma_{Nm}^{mech}}\vc t_{Nm}\cdot(\zz_m-\uu_m)
+\int_{\Omega_f}\delta\ff_f\cdot(\zz_f-\uu_f) + \int_{\Gamma_{Nf}^{mech}}\delta\vc t_{Nf}\cdot(\zz_f-\uu_f).
}\end{linenomath}
Splitting the total stress $\sigmapor_*$, $*\in\{m,f\}$, into the mechanical and hydraulic part, we obtain the variational inequality \cref{eq:weak_var_ineq}.

\paragraph{Weak form of conservation of mass \cref{eq:weak_flow}} 
The weak form of \cref{eq:fl,eq:mixed_dim_problem_fl_frac} is obtained directly by using test functions $q=(q_m,q_f)\in \mathcal Q$:
\begin{linenomath}\begin{align*}
\int_{\Omega_m}\left(\partial_t(S_m p_m + \alpha_m\varrho g\div\uu_m) + \div\vv_m\right)q_m &= \int_{\Omega_m}s_m q_m\\
\int_{\Omega_f}\left(\partial_t\left(\delta S_f p_f + \delta\alpha_f\varrho g\fdiv\uu\right) + \fdiv(\delta\vv_m,\vv_f)\right)q_f &= \int_{\Omega_f}\delta s_f q_f.
\end{align*}\end{linenomath}
This can be written in the compact form \cref{eq:weak_flow}.

\paragraph{Weak form of the Darcy law \cref{eq:weak_darcy}}
Let $\ww=(\ww_m,\ww_f)\in \mathcal W(0)$.
Then the Darcy law in $\Omega_m$, Green's theorem and the boundary conditions on $\Gamma_{Dm}^{flow}$ and $\Gamma_{Nm}^{flow}$ imply:
\begin{linenomath}\ml{\label{eq:darcy_weak_m}
\int_{\Omega_m}\tn K_m^{-1}\vv_m\cdot\ww_m 
=  - \int_{\Omega_m}(\nabla p_m+\vc g)\cdot\ww_m\\
=  - \int_{\Omega_m}\vc g\cdot\ww_m -\int_{\Gamma_{Dm}^{flow}}p_{Dm}\nn\cdot\ww_m + \int_{\Omega_m}p_m\div\ww_m + \int_{\Omega_f}2\avg{p_m\ww_m\cdot\nnu}.
}\end{linenomath}
Since the fracture conductivity tensor $\tn K_f$ is isotropic, we can express the flux continuity \cref{eq:cont_flux} as follows: 
\begin{linenomath}\eqs{ \vv_m^\pm\cdot\nnu^\pm = -\frac{2k_f(\uu)}\delta(p_m^\pm-p_f) - k_f(\uu)\vc g\cdot\nnu^\pm. }\end{linenomath}
Using this relation, we rewrite the last term from \cref{eq:darcy_weak_m}:
\begin{linenomath}\mls{
\int_{\Omega_f}2\avg{p_m\ww_m\cdot\nnu}
= \int_{\Omega_f} 2\avg{ (p_m-p_f)\ww_m\cdot\nnu } + \int_{\Omega_f} 2 p_f\avg{\ww_m\cdot\nnu}\\
= -\int_{\Omega_f}\frac\delta{k_f}\avg{(\vv_m\cdot\nnu) (\ww_m\cdot\nnu)}
-\int_{\Omega_f}\delta\avg{(\vc g\cdot\nnu)(\ww_m\cdot\nnu)}\\
+ \int_{\Omega_f} 2p_f\avg{\ww_m\cdot\nnu}.
}\end{linenomath}

Next, from the Darcy law in $\Omega_f$ and the appropriate boundary conditions, we get:
\begin{linenomath}\mls{
\int_{\Omega_f}\frac1\delta\tn K_f^{-1}(\uu)\vv_f\cdot\ww_f
= - \int_{\Omega_f}(\nabla_\tau p_f + \vc g_\tau)\cdot\ww_f\\
= -\int_{\Omega_f}\vc g_\tau\cdot\ww_f - \int_{\Gamma_{Df}^{flow}} p_{Df}\nn\cdot\ww_f + \int_{\Omega_f}p_f\div_\tau\ww_f.
}\end{linenomath}
Altogether we obtain:
\begin{linenomath}\mls{
\int_{\Omega_m}\tn K_m^{-1}\vv_m\cdot\ww_m + \int_{\Omega_f}\frac1\delta k_f^{-1}(\uu)\vv_f\cdot\ww_f
+ \int_{\Omega_f}\frac\delta{k_f}\avg{(\vv_m\cdot\nnu)(\ww_m\cdot\nnu)}\\
- \int_{\Omega_m}p_m\div\ww_m
- \int_{\Omega_f} p_f\fdiv(\delta\ww_m,\ww_f)\\
= - \int_{\Omega_m}\vc g\cdot\ww_m - \int_{\Omega_f}\vc g_\tau\cdot\ww_f
-\int_{\Omega_f}\delta\avg{(\vc g\cdot\nnu)(\ww_m\cdot\nnu)}\\
- \int_{\Gamma_{Dm}^{flow}}p_{Dm}\nn\cdot\ww_m - \int_{\Gamma_{Df}^{flow}}p_{Df}\nn\cdot\ww_f.
}\end{linenomath}
The last equation can be rewritten as \cref{eq:weak_darcy}.

\section{List of symbols and abbreviations}

Symbols with subscript ``m'' and ``f'' denote quantities related to the rock matrix and fractures, respectively.

\noindent
\begin{longtable}[l]{ll}
\multicolumn{2}{l}{\bf Latin symbols}\\
\hline & \\[-1.5ex]
$a_f$ & fracture aperture\\
$\tn C_m$, $\tn C_f$ & elasticity tensor\\
$E_m$, $E_f$ & Young's modulus\\
$\ff_m$, $\ff_f$ & load\\
$g$ & gravitational acceleration\\
$\vc g$ & gradient of vertical distance\\
$H^1$ & Sobolev space\\
$\Hdiv$ & space of vector-valued functions with square-integrable value and divergence\\
$I$ & time interval\\
$\tn I$ & second order unit tensor\\
$\tn K_m$, $\tn K_f$ & hydraulic conductivity (tensor)\\
$k_m$, $k_f$ & hydraulic conductivity (scalar)\\
$L^2$ & Lebesgue space\\
$\nn$ & unit outward normal vector\\
$p$, $p_m$, $p_f$ & pressure head\\
$p_0$, $p_{0m}$, $p_{0f}$ & initial pressure head\\
$p_D$, $p_{Dm}$, $p_{Df}$ & pressure head on the Dirichlet boundary\\
$P_k$ & space of polynomials with degree up to $k$\\
$RT_k$ & $k$-th order Raviart-Thomas space\\
$\Real$ & set of all real numbers\\
$S_m$, $S_f$ & storativity\\
$\vc t_N$, $\vc t_{Nm}$, $\vc t_{Nf}$ & traction on the Neumann boundary\\
$\uu$, $\uu_m$, $\uu_f$ & displacement\\
$\uu_D$, $\uu_{Dm}$, $\uu_{Df}$ & displacement on the Dirichlet boundary\\
$\vv_m$, $\vv_f$ & Darcian flux (velocity) in the rock matrix, Darcian flux through fractures\\
$v_N$, $v_{Nm}$, $v_{Nf}$ & flux on the Neumann boundary\\
\end{longtable}

\bigskip

\noindent
\begin{longtable}[l]{ll}
\multicolumn{2}{l}{\bf Greek symbols}\\
\hline & \\[-1.5ex]
$\alpha_m$, $\alpha_f$ & Biot-Willis coefficient\\
$\Gamma_{D\{m,f\}}^{\{flow,mech\}}$ & Dirichlet boundary for flow/mechanics\\
$\Gamma_{N\{m,f\}}^{\{flow,mech\}}$ & Neumann boundary for flow/mechanics\\
$\delta$, $\delta_{min}$ & fracture width, minimal fracture width\\
$\eps$ & symmetric gradient\\
$\eta$ & fracture roughness\\
$\vc\Lambda^\pm$ & contact force\\
$\mu$ & fluid viscosity\\
$\nnu$ & unit normal vector to fracture plane\\
$\nu_m$, $\nu_f$ & Poisson's ratio\\
$\Omega$, $\Omega_m$, $\Omega_f$ & computational domain, domain for rock matrix, domain for fractures\\
$\varrho$ & fluid density\\
$\vc\sigma_m$, $\vc\sigma_f$ & effective stress\\
$\vc\sigma_{0m}$, $\vc\sigma_{0f}$ & initial stress\\
$\sigmapor$, $\sigmapor_m$, $\sigmapor_f$ & total stress\\
\end{longtable}

\bigskip

\noindent
\begin{longtable}[l]{ll}
\multicolumn{2}{l}{\bf Other symbols}\\
\hline & \\[-1.5ex]
$\nabla_\tau$ & tangential gradient\\
$\fgrad$, $\fgrad^\pm$ & semi-discrete gradient\\
$\div_\tau$ & tangential divergence\\
$\fdiv$ & semi-discrete divergence\\
$\feps$, $\feps^\pm$ & semi-discrete symmetric gradient\\
$\avg{\ }$ & average of a quantity over fracture sides\\
$(\ )^+$, $(\ )^-$ & value on positive/negative fracture side\\
\end{longtable}

\bigskip

\noindent
\begin{longtable}[l]{ll}
\multicolumn{2}{l}{\bf Abbreviations}\\
\hline & \\[-1.5ex]
DFM & discrete fracture-matrix (model)\\
DFN & discrete fracture network\\
FETI & finite element tearing and interconnecting (method)\\
GFEM & generalized finite element method\\
MPGP & modified proportioning with gradient projections (method)\\
QP & quadratic programming\\
TSX & tunnel sealing experiment\\
XFEM & extended finite element method\\
\end{longtable}

\bibliographystyle{elsarticle-num}

\bibliography{ref}

\end{document}